\documentclass[10pt,leqno]{amsart}
\usepackage{amsfonts,amssymb}
\usepackage{tikz, subfigure, xcolor} 

\usepackage{amsmath}
\usepackage{amsthm}
\usepackage[alphabetic]{amsrefs}
\usepackage{qsymbols}
\usepackage{latexsym}
\usepackage{chngcntr}
\usepackage[noadjust]{cite}
\usepackage{paralist}
\usepackage{doi}
\usepackage{comment}
\usepackage{dsfont}
\usepackage{esint}

\usepackage{color}

\providecommand{\norm}[1]{\left\| #1 \right\|} 

\newtheorem{theorem}{Theorem}[section]

\newtheorem{lemma}[theorem]{Lemma}
\newtheorem{proposition}[theorem]{Proposition}
\newtheorem{corollary}[theorem]{Corollary}
\theoremstyle{definition}

\newtheorem{remark}[theorem]{Remark}
\newtheorem{remarks}[theorem]{Remarks}

\theoremstyle{plain}

\renewcommand{\d}{\operatorname{d}\!} 
\DeclareMathOperator{\dt}{\d t}		
\DeclareMathOperator{\dW}{\d W} 	

\newcommand{\IR}{\mathbb{R}}
\newcommand{\IC}{\mathbb{C}}
\newcommand{\IN}{\mathbb{N}}
\newcommand{\IZ}{\mathbb{Z}}

\newcommand{\IP}{\mathbb{P}}

\newcommand{\E}{\mathbb{E}}

\newcommand{\cH}{\mathcal{H}}

\newcommand{\cF}{\mathcal{F}}

\newcommand{\cD}{\mathcal{D}}

\newcommand{\abs}[1]{\left|#1\right|}

\newcommand\restr[2]{{
		\left.\kern-\nulldelimiterspace 
		#1 
		\vphantom{\big|} 
		\right|_{#2} 
}}

\renewcommand{\d}{\mathrm{d}}

\renewcommand{\div}{\mathrm{div}\,}    

\newcommand{\lpso}{L^p_{\overline{\sigma}}(\cD;\cH)}
\newcommand{\ltso}{L^2_{\overline{\sigma}}(\cD;\cH)}
\newcommand{\Ee}{\mathbb{E}}

\numberwithin{equation}{section} 

\title[The primitive equations with stochastic wind driven boundary conditions]{The primitive equations with stochastic wind driven boundary conditions}



\author[Binz]{Tim Binz} 
\address{Department of Mathematics,
	TU Darmstadt, Schlossgartenstr. 7, 64289 Darmstadt, Germany}
\email{binz@mathematik.tu-darmstadt.de}
\email{hieber@mathematik.tu-darmstadt.de}
\email{msaal@mathematik.tu-darmstadt.de}

\author[Hieber]{Matthias Hieber} 

\author[Hussein]{Amru Hussein}
\address{Department of Mathematics,
	TU Kaiserslautern, Paul-Ehrlich-Stra{\ss}e 31,
	67663 Kaiserslautern, Germany}
\email{hussein@mathematik.uni-kl.de}

\author[Saal]{Martin Saal}



\begin{document}

\begin{abstract}
	The primitive equations for geophysical flows are studied under the influence of {\em stochastic wind driven boundary conditions} modeled by a cylindrical Wiener process. 
	We adapt an approach by  Da Prato and Zabczyk for stochastic boundary value problems to define a notion of solutions.
	Then a  rigorous treatment of these stochastic boundary conditions, which combines stochastic and deterministic methods,  yields that these equations admit a unique, local pathwise solution within the 
	anisotropic $L^q_t$-$H^{-1,p}_zL^p_{xy}$-setting. 
	This solution is constructed in critical spaces.    
\end{abstract}
\maketitle


\section{Introduction}
Consider the primitive  equations in a cylindrical domain $\cD = G \times (-h,0) \subset \IR^3$,  where 
$G=(0,1)\times(0,1)$ and  $h>0$. Let us denote by $v\colon\cD \times (0,T)\rightarrow \IR^2$ the horizontal velocity of the fluid 
and by $p_s\colon G \times (0,T) \rightarrow \IR$ its surface pressure on a time interval $(0,T)$, where $T>0$. 
We consider the set of equations  
\begin{align}\label{eq:primequiv}
\left\{
\begin{array}{rll}
\partial_t v + v \cdot \nabla_H v + w(v) \cdot \partial_z v - \Delta v + \nabla_H p_s  & = f, &\text{ in } \cD \times (0,T),  \\
\mathrm{div}_H \overline{v} & = 0, &\text{ in } \cD \times (0,T), \\
v(0) & = v_0, &\text{ in } \cD,
\end{array}\right.
\end{align}
where  $\overline v(x,y) = \frac1h \int_{-h}^0 v(x,y, \xi)d\xi$, and  the vertical velocity $w=w(v)$ with $w(x,y,-h)=w(x,y,0)=0$ is given by $w(v)(x,y,z) = -\int_{-h}^z \mathrm{div}_H v(x,y, \xi)d\xi$. Here  $(x,y)\in G$ denote the horizontal coordinates and 
$z\in (-h,0)$ the vertical one. 
There exist  several equivalent 
formulations of the primitive equations, depending on whether the vertical velocity $w=w(v)$ is completely substituted by   the horizontal velocity $v$ and the full pressure 
by the surface pressure, respectively, compare e.g. \cite{HK16}. 

It is the  aim of this article to study the above set of equations subject to {\em stochastic wind driven boundary conditions}. These boundary conditions on the atmosphere-ocean interface describe 
the balance of the shear stress of the ocean and the horizontal wind force. In contrast to previous works on deterministic wind driven boundary conditions described by Lions, Temam and Wang in 
\cite{Lionsetal1993}, we investigate here for the first time stochastic wind driven boundary conditions of the form $\partial_z v = h_b \partial \omega$ for a given cylindrical Wiener process
$\omega(t) = \sum_{n=1^\infty} <g,e_n> W_b(t)e_n$. For a precise definition of this condition, we refer to Section~\ref{sec:prel_stoch} below.          

Wind driven boundary conditions for the coupled atmosphere and ocean  primitive equations within the deterministic setting were introduced and studied by Lions, Temam and Wang in their fundamental article 
\cite{Lionsetal1993}. For related  results concerning deterministic wind driven boundary conditions for the Navier-Stokes equations we refer to the work of Desjardins and Grenier \cite{DG00}, Bresch and 
Simon \cite{BS01} and Dalibard and Saint-Raymond \cite{DS09}.
Here, the equations \eqref{eq:primequiv} are supplemented by mixed boundary conditions on $\Gamma_{u} = G \times \{0\}$, $\Gamma_b = G \times \{-h\}$ and $\Gamma_l = \partial G \times (-h,0)$
of the form 
\begin{align}
&v, p_s  \hbox{ are periodic } \hbox{on } \Gamma_l \times (0,T),  \label{bcl}\\
\partial_z v =& 0 \hbox{ on } \Gamma_b \times (0,T),  \label{bcb}\\
\partial_z v =& c\varrho^{air}(v^{air}-v)\cdot |v^{air}-v|  \hbox{ on } \Gamma_u \times (0,T). \label{bcu}
\end{align}
Here $v^{air}$ denotes the velocity of the wind, $\varrho^{air}$ the density of the atmosphere and $c$ the drag coefficient. The boundary condition \eqref{bcu} is interpreted as the physical law 
describing the driving mechanism on the atmosphere-ocean interface as a balance of the shear stress of the ocean and the horizontal wind force. Indeed, the shear stress of the ocean, i.e. 
the tangential component of the stress tensor is given by $\partial_z v + \nabla_H w$, which due to the flatness of the interface, i.e. $w=0$ on $\Gamma_u$, equals $\partial_z v$, for details see 
\cite{Lionsetal1993}. 

At first glance a natural boundary  condition on the interface would be an adherence condition, i.e. $v=v^{air}$, at the interface. These conditions are, however, not being used due to the occurrence 
of boundary layers in the atmosphere and in the ocean at the surface. The above condition \eqref{bcu} takes into account these boundary layers.  Since the velocity of air is much slower than the one of the 
ocean, the term $v$ is frequently neglected and the condition 
\begin{equation}\label{bcusimple}
\partial_z v = c\varrho^{air}v^{air}\cdot |v^{air}|  \hbox{ on } \Gamma_u \times (0,T)
\end{equation}
is used instead, see e.g. \cite{TZ96,GM97}. 


In this article we extend the above setting in the following  direction: we introduce  wind driven {\em stochastic} boundary conditions on the surface of the ocean and analyze the 
primitive equation subject to these boundary conditions as an SPDE. Stochastic wind driven boundary conditions have been considered before within the setting of  the shallow water equations e.g. 
by Cessi and Louazei \cite{CL01} from a modeling point of view. For numerical results and statistical analysis of wind stress time series in the context of the Ekman equation, we refer 
to the work \cite{BCP15} of Buffoni, Cappeletti and Picco. Our result seems to be the first rigorous result concerning stochastic boundary conditions driven by wind.  A similar setting has been considered very recently for the 2-dimensional Navier-Stokes equations in \cite{agresti2023global}.
%
More precisely, given a cylindrical Wiener process $W$ on a separable Hilbert space $\cH$ with respect to a filtration $\mathcal{F}$ and adapted functions $H_f$ and $h_b$, we consider for the horizontal velocity of the fluid
$V\colon\Omega \times \cD\times (0,T) \rightarrow \IR^2$  
and the surface pressure $P_s\colon \Omega \times G \times (0,T)\rightarrow \IR$, where $(\Omega, \mathcal{A},P)$ is a probability space endowed with the filtration $\mathcal{F}$, 
the equations 
\begin{align}\label{eq:pe_stochi}
\left\{
\begin{array}{rll}
\d V + (V \cdot \nabla_H V + w(V) \cdot \partial_z V - \Delta V + \nabla_H P_s) \d t & = H_f\d W, &\text{ in } \cD \times (0,T),  \\
\mathrm{div}_H \overline{V} & = 0, &\text{ in } \cD \times (0,T), \\
V(0) & = V_0, &\text{ in } \cD, 
\end{array}\right.
\end{align}
subject to boundary conditions \eqref{bcl} and \eqref{bcb}, but where the deterministic condition \eqref{bcu} or \eqref{bcusimple} is replaced by  a stochastic boundary condition modeling the wind as 
\begin{align}\label{bc:stochi}
\partial_z V = h_b \, \partial_t \omega \quad \hbox{ on } \Gamma_{u} \times (0,T).
\end{align} 
Here $h_b$ is a function defined on $\Gamma_u\times (0,T)$ and we assume that  $\omega$ can be written as
\begin{align}\label{defomega}
\omega(t) = \sum_{n=1}^{\infty} <g,e_n>  W_b(t) e_n, 
\end{align}
where $g$ is a suitable function defined on $\Gamma_u$, $W_b$ is  another cylindrical Wiener process on $\cH$ with respect to the filtration $\mathcal{F}$, and $(e_n)$ is an orthonormal basis of 
$\cH$.   

Our strategy to  prove the existence of a unique, local pathwise solution for equations \eqref{eq:pe_stochi} and \eqref{bc:stochi} is based on a combination of stochastic and deterministic maximal $L^q$-regularity  methods.
It can be summarized as follows: First, in order to eliminate the pressure term we apply the hydrostatic Helmholtz projection $\mathbb{P}$ to equation \eqref{eq:pe_stochi} (compare e.g. \cite{HK16}) and 
rewrite the stochastic primitive equations as a semilinear stochastic evolution equation 
of the form   
\begin{align}\label{eq:primaddnoise}
&\d V+ A V \dt =F(V,V) \dt + H_f \d W, \quad V(0)=V_0.
\end{align}
Here  $A$ denotes the hydrostatic Stokes operator defined  as $A=-\mathbb{P}\Delta$ and 
$F(\cdot,\cdot)$ is the bilinear convection term. In \cite{HK16} this has been done in the space $L^p_{\overline{\sigma}}(\cD)$. 
Here, however, the choice of the ground space and correspondingly the domain $D(A)$ will give us some freedom, compare Section~\ref{sec:prel} below,  needed to include the stochastic boundary conditions.   

Secondly, we rewrite the stochastic boundary condition as a forcing term. Indeed,
a general obstacle to handle stochastic  boundary conditions is a lack of regularity to make sense of \eqref{bc:stochi}. Here, 
we adapt  an approach due to Da Prato and Zabczyk \cite{DZ96} for stochastic boundary conditions  to the given situation:
A solution $V$
to equation \eqref{eq:pe_stochi}
subject to \eqref{bcl}, \eqref{bcb} and the stochastic condition \eqref{bc:stochi} is expressed by a solution to the equation
\begin{align}
\d Z_b(t) + A Z_b(t) \dt = A[\mathcal{N} h_b(t) g] \dW_b(t), \quad Z_b(0)=0,
\end{align}
subject to the  boundary conditions  
\begin{align}\label{eq:Zb_bc}
Z_b \hbox{ are periodic } \hbox{on } \Gamma_l \times (0,T),  \quad  \hbox{and}\quad
\partial_z Z_b = 0 \hbox{ on } \Gamma_{u} \cup \Gamma_{b} \times (0,T).
\end{align}
Here $\mathcal{N}$ denotes the so-called Neumann operator mapping deterministic inhomogeneous boundary data to the solution of the  associated stationary hydrostatic Stokes problem. This hydrostatic Neumann operator is constructed in Subsection~\ref{subsec:Neumann} below. This construction allows us to view the 
stochastic boundary condition as a stochastic forcing term. For a similar approach within the setting of  parabolic equations in divergence form we also refer  to \cite{SchnaubeltVeraar2011}.
Here, we call $V$ a solution to \eqref{eq:pe_stochi}
subject to \eqref{bcl}, \eqref{bcb} and the stochastic condition \eqref{bc:stochi} 
if $V_b:=V-Z_b$ solves the equation
\begin{align}\label{eq:Vb}
\d V_b+ A V_b \dt =F(V_b+Z_b,V_b+Z_b) \dt + H_f \d W, \quad V(0)=V_0
\end{align}
subject to the same homogeneous boundary conditions \eqref{eq:Zb_bc}.  
As discussed in detail by Da Prato and Zabczyk in~\cite[Section 13]{DZ96}, if this solution $V$  were sufficiently regular, then it would constitute a solution to the actual  inhomogeneous problem.
We note that $Z_b$ is given by
\begin{align}\label{eq:Zb_intro}
Z_b(t):=A\int_0^t e^{-(t-s)A} [\mathcal{N} h_b(t)g] \dW_b(t), \quad t>0.
\end{align}

Thirdly, we  investigate the solution $Z_f$ of the linearized system with linear noise 
\begin{align*}
&\d Z_f+ A_p Z_f \dt = H_f \d W, \quad Z_f(0)=Z_0.
\end{align*}
Subsequently, we consider  pathwise the remainder term $v:=V-Z$ for $Z:=Z_f + Z_b$, which solves (almost surely) the {\em deterministic and nonautonomous},  semilinear evolution equation
\begin{align}\label{eq:primpath}
&\partial_t v+ A v =F(v+Z,v+Z), \quad v(0)=v_0,
\end{align}
where the initial value $Z_0$ contains the probabilistic part of the initial value $V_0$ and  $v_0:=V_0-Z_0$ its 
deterministic part.  
This equation will 
be treated by the theory of semilinear evolution equations in critical spaces due to  
Pr\"uss, Simonett and Wilke \cite{PruessSimonettWilke, PruessWilke2016}. 
This will enables us 
to  prove the existence of a {\em unique, local} solution to \eqref{eq:primpath}, 
whereby we use the theory of time weighted maximal $L^q_t$-regularity.

After these reformulations, the main challenge is now to assure sufficient
regularity properties of $Z$ which allow us to  solve \eqref{eq:primpath} pathwise in the deterministic maximal regularity framework.
To this end, we use results on maximal stochastic regularity due to van Neerven, Veraar and Weis \cite{NeervenVeraarWeis}. The latter are applicable due to the fact 
that $A$ admits a bounded $H^\infty$-calculus in $L^p_{\overline{\sigma}}(\cD)$, see \cite{GGHHK17}, which carries over to further anisotropic spaces by isomorphy. 
By the theory of stochastic maximal integrals the regularity stochastic convolutions of the form \eqref{eq:Zb_intro} increases by an order of $A^{1/2}$ provided that $A\mathcal{N}(h_b(\cdot)g)$ lies in the corresponding ground space. 
Now,
\begin{align*}
\mathcal{N}(\cdot)\colon W^{1-1/r+m,p}(\Gamma_u;\cH)^2\rightarrow H^{2+m,p}(\cD;\cH)^2 = H^{2+m,p}_zL^p_{xy} \cap L^{p}_zH^{2+m,p}_{xy},
\end{align*}
where it is essential to observe that here $[\mathcal{N} h_b(t)g]$ satisfies the inhomogeneous boundary conditions
\eqref{eq:bc_l_and_u}. Hence, $[\mathcal{N} h_b(t)g]\in D(A)$ holds
only if we shift the linear operator to a weaker setting. 
Our strategy here, is to make this shift only with respect to  the vertical $z$-directions, where 
the inhomogeneous boundary condition is rendered incognisable in $H^{s,p}_z$ below the threshold of  $s<1+1/p$. Heuristically, this gives a vertical regularity for $Z_b$ of order $H^{s-1,p}_z$ with $s-1<1/p$, and thus we have less than one vertical derivative available. This makes it necessary to consider anisotropic ground spaces of the type $H^{s-2,p}_{z}H^{m,p}_{xy}$ for some $s$ as above and $m\geq 0$.

For the linear part of the equations shifting the scales is possible by abstract extrapolation scales. However, the limiting factor is the availability of the corresponding  non-linear estimates in such weak settings. In the Navier-Stokes equations  the non-linearity can take the form $\div (u \otimes u)$ thus allowing for a shift to spaces such as $H^{-1,p}(\cD)$. In contrast, for the primitive equations the non-linearity takes the form $\div ((v,w(v)) \otimes v)$, where $w(v)$ still contains derivatives of $v$ by \eqref{eq:w} below. However, we still achieve non-linear estimates in the anisotropically weak space $H^{-1,p}_{z}L^{p}_{xy}$ which  lead to unique strong solutions in this setting. This should not be confused with weak solutions 
or the so-called $z$-weak solutions for the primitive equations, compare e.g. \cite{Ju2017}. 
In the estimate on  $F(Z_b,Z_b)$, where $Z_b$ is as in \eqref{eq:Zb_intro} and has less than one derivative in $z$-direction, we can compensate for this by adding some regularity in horizontal direction by some $m>0$.
Thus, we have to assume that $h_b(\cdot)g\in L_\sigma^{2q}(0,T; W^{1+1/p+\varepsilon}(\Gamma_u;\cH)^2)$ to solve \eqref{eq:primpath} in $L_{\mu}^{q}(0,T;H_{z}^
{-1,p}L_{xy}^p)$ for suitable time-weights $\sigma$ and $\mu$, and for $p\in (3,4]$ and for some $q$.

To extend the local solutions to a global one suitable \textit{a priori} bounds are needed. In the deterministic setting \textit{a priori} bounds in the maximal $L^2_t$-$L_x^2$-setting are available by the classical result of Cao and Titi, cf. \cite{CaoTiti2007}. However, for the setting needed here, {a priori} bounds in $H^{-1,p}_zL^p_{xy}$-spaces would be needed which are -- as of now --  not available, 
and which are subject of a future work.

Recently,  Agresti and Veraar \cite{AgrestiVeraar, AgrestiVeraarII} developed a {\em local} theory of critical spaces for stochastic evolution equations analogously to the ideas in \cite{PruessSimonettWilke} 
for deterministic equations. One  major difference is that due to the 
weaker smoothing of the stochastic convolution, the conditions 
on the weights used by them are more restrictive as in the deterministic case. By  using our approach when solving \eqref{eq:primpath} we are able to allow  spatially rougher data $v_0$ than one 
could handle by considering \eqref{eq:primpath} in the context of stochastic critical spaces.


Recall that the mathematical analysis of the {\em deterministic} primitive equations has been pioneered by Lions, Teman and Wang  in their articles \cite{Lionsetal1992, Lionsetal1992_b, Lionsetal1993}, where 
the existence of a global, weak solution to the primitive equations is proven. For global weak solutions  subject to \eqref{bcu} we refer to  \cite{PTZ09}. 
The uniqueness property of weak solutions for initial data in $L^2$ remains an open problem until today. 
A landmark result on the global strong well-posedness of the {\em deterministic} primitive equations subject to homogeneous Neumann conditions for initial data in $H^1$ was shown by Cao and Titi in 
\cite{CaoTiti2007} by the method of energy estimates. For mixed Dirchlet-Neumann conditions we refer to the work  \cite{KZ07} of Kukavica and Ziane.    
A different approach to the deterministic  primitive equations, based on methods of evolution equations, has been introduced  in \cite{GGHHK20b, HK16}. This approach is based on the  
hydrostatic Stokes operator $A$ and the hydrostatic Stokes semigroup defined  for $p \in (1,\infty)$ on the hydrostatic solenoidal $L^p$-spaces. 
For a survey on results concerning the deterministic primitive equations using the approach of  energy estimates, we refer to \cite{LiTiti2016}; for a survey concerning the approach based on 
evolution equations, see \cite{HH20}.

The three dimensional \emph{stochastic} primitive equations with deterministic boundary conditions but stochastic forcing term have been studied before by several authors. 
Indeed, for the situation of  additive noise, there are existence and uniqueness results for pathwise, strong 
solutions within the $L^2$-setting;  see \cite{GuoHuang2009}. They consider deterministic initial data in $H^1(\cD)$ and choose Neumann boundary conditions on the bottom and top.

A global well-posedness result for pathwise strong solutions of the primitive equations with deterministic and homogeneous boundary conditions was established for multiplicative 
white noise in time in  \cite{DebusscheGlattholtzTemam2011, DebusscheGlattholtzTemamZiane2012} and later under weaker assumptions on the noise 
in \cite{BS20}. Here  Neumann boundary conditions are used  for the top and the bottom in \cite{DebusscheGlattholtzTemamZiane2012} and a Dirichlet boundary
condition on the bottom combined with a mixed Dirichlet-Neumann boundary conditions on the top in \cite{DebusscheGlattholtzTemam2011}. 
Further results concerning  the existence of ergodic invariant measures, weak-martingale solutions and Markov selection  
were shown in \cite{DongZhang2017} and  \cite{GlattholtzKukavicaVicolZiane2014}. For results in  two dimensions,  see e.g. \cite{GlattholtzTemam2011}.
It seems that here for the first time stochastic boundary conditions are studied for the primitive equations.

In the following, we elaborate on our strategy outlined above: We fix the setting and some notation in Section~\ref{sec:pe}.
In order to adapt the approach by Da Prato and Zabczyk  to the stochastic primitive equations with stochastic boundary conditions 
we have to analyse the linearized equations carefully  for both the deterministic setting in Section~\ref{sec:prel} and then the stochastic setting in Section~\ref{sec:prel_stoch}.
Then we we are in the position to define our notion of solution and to  give our main  result on the local existence and uniqueness of solutions and its proof in Section~\ref{sec:main}. 






\section{First formulation of the stochastic primitive equations}\label{sec:pe}
We consider  the stochastic primitive equations in the isothermal setting on a cylindrical spatial domain 
\begin{align*}
\cD =  G \times (-h,0) \subset \IR^3\quad \hbox{with}\quad  G=(0,1)\times(0,1), \quad \hbox{where }h>0,
\end{align*}  
on a time interval $(0,T)$ with $T>0$, and a probability space  $(\Omega, \mathcal{A},P)$. The upper, bottom and lateral parts of the boundary $\partial\cD$, respectively, are denoted by
\begin{align*}
\Gamma_{u} = G \times \{0\}, \quad \Gamma_b = G \times \{-h\}, \quad \hbox{and} \quad \Gamma_l = \partial G \times (-h,0).
\end{align*}
The unknowns are the horizontal velocity of the fluid and the surface pressure 
\begin{align*}
V\colon\Omega\times \cD \times (0,T) \rightarrow \IR^2 \quad \hbox{and} \quad P_s\colon \Omega\times \cD \times (0,T) \rightarrow \IR,
\end{align*}
respectively. 
These are governed by the following 
already reformulated stochastic primitive equations
\begin{align}\label{eq:pe_stoch}
\left\{
\begin{array}{rll}
\d V + (V \cdot \nabla_H V + w(V) \cdot \partial_z V - \Delta V + \nabla_H P_s) \d t  
&   = H_f \, \d W, &\text{ in } \cD \times (0,T),  \\
\mathrm{div}_H \overline{V} & = 0, &\text{ in } \cD \times (0,T), \\
V(0) & = V_0, &\text{ in } \cD 
\end{array}\right.
\end{align}
for given initial data $V_0\colon \Omega \times  \cD \rightarrow \IR^2$ and  stochastic forcing $H_f \, \d W$  defined by a cylindrical Wiener process and a given function $H_f$ made precise in Subsection~\ref{subsec:Wiener} below. 
Here $(x,y)\in G$ denote the horizontal coordinates, $z\in (-h,0)$ the vertical one, and 
\begin{align*}
\Delta = \partial_x^2 + \partial_y^2+ \partial_z^2, \quad \nabla_H  = (\partial_x, \partial_y)^T, \quad \mathrm{div}_H V= \partial_x V_1+ \partial_y V_2, \\ \hbox{and} \quad \overline{V}:=\frac{1}{h}\int_{-h}^0 V(\cdot,\cdot, \xi)d\xi.
\end{align*}
The vertical velocity $w=w(V)$ in \eqref{eq:pe_stoch} is given by
\begin{align}\label{eq:w}
w(V)(x,y,z) = -\int_{-h}^z \mathrm{div}_H V(x,y, \xi)d\xi, 
\end{align}
because of the boundary condition 
\begin{align*}
w(V) &= 0 \hbox{ on } \Gamma_{b}\cup \Gamma_u\times (0,T).
\end{align*}
The equations \eqref{eq:pe_stoch} are supplemented by the boundary conditions on the lateral part
\begin{align}\label{eq:bc_l}
V, P_s  \hbox{ are periodic } \hbox{on } \Gamma_l\times (0,T),
\end{align}
on the bottom part 
\begin{align}\label{eq:bc_Neumann}
\partial_z V &= 0 \hbox{ on } \Gamma_{b}\times (0,T),  
\end{align}
and a stochastic forcing term modeling the wind is imposed on the upper part by 
\begin{align}\label{eq:bcstoch}
\partial_z V = h_b \partial_t \omega \hbox{ on } \Gamma_{u} \times (0,T).
\end{align}
Here, $h_b\colon \Gamma_u \times (0,T) \rightarrow \IR$ is a given real valued function, the assumptions on which are made precise in Subsection~\ref{subsec:boundary} below, 
and $\partial_t \omega$ stands for a noise term defined by a cylindrical Wiener process
which is  discussed in Subsection~\ref{subsec:Wiener} below.

\section{Linear deterministic theory}\label{sec:prel}


\subsection{Anisotropic  $H_{z}^{s,p}H^{m,p}_{xy}$-spaces}

General vector valued function space and distributions in anisotropic spaces are discussed systematically by Amann in \cite{Amann2019} for corner-domains in $\IR^n$, and the anisotropic scalar case can be found in \cite[Section 10.1]{Triebel}. In both instances it is assumed that the differentiability 
has a fixed sign for all directions. Here, we aim to extend this to Sobolev spaces   with negative order in $z$-direction and positive order in $x$-$y$-direction.
Note that the vector valued periodic case in isotropic  spaces for positive differentiability is discussed also for instance by Arendt and Bu in \cite{AB:02} and by Nau in \cite{Nau2012}. 

For a separable Hilbert space $\cH$, $p\in(1,\infty)$ and $s\in [0,\infty)$, we define the periodic isotropic $\cH$-valued Bessel potential spaces with respect to the horizontal variables dealing with negative orders via duality
\begin{align*}
H^{s,p}_{per}(G;\cH) := \overline{C^{\infty}_{per}(\overline{G};\cH)}^{\norm{\cdot}_{H^{s,p}(G;\cH)}} \quad \hbox{and} \quad H^{-s,p}_{per}(G;\cH):=(H^{s,p'}_{per}(G;\cH))'.
\end{align*}
Here, $1/p+1/p'=1$, $C^{\infty}_{per}(\overline{G};\cH)$ 
is the space of smooth  on $\overline{G}$ with values in $\cH$ which are periodic of any order, cf. also \cite[Chapter 2]{Nau2012}, and we have identified $\cH$ and it dual $\cH'$ via the Riesz isomorphism. For $s=0$ we have $H^{0,p}= L^p$. 

To define anisotropic function spaces, consider first more generally
for a given Banach space $X$ the space of test functions
\begin{align*}
\cD_{per,ev}(X):=
\{v\in C^{\infty}([-h,h];X)\colon v \hbox{ periodic and even w.r.t. the $z$-direction}\},
\end{align*}
and we set $C_{per, ev}^{\infty}([-h,h];X):=\cD_{per,ev}(X)$.
The image of the restriction operator
\begin{align*}
R \colon \cD_{per,ev}(X) \rightarrow \cD_{N}(X):= \cD_{per,ev}(X)\vert_{(-h,0)}, \quad v\mapsto v\vert_{(-h,0)}
\end{align*}
satisfies Neumann-type boundary conditions on top and bottom, that is,
\begin{align*}
\partial_z^{2n-1}v(-h)=\partial_z^{2n-1}v(0)=0\quad \hbox{for }v\in \cD_{N}(X) \hbox{  and } n\in \IN,
\end{align*}
compare also \cite[Chapter 7]{Nau2012}.
Moreover, the even extension operator 
\begin{align*}
E_{ev}\colon \cD_{N}(X) \rightarrow \cD_{per,ev}(X), \quad E_{ev}v(z):=\begin{cases}
v(z), & z\in [-h,0], \\
v(-z), & z\in [0,h],
\end{cases}
\end{align*}
is well-defined. The operators $E_{ev}$ and $R$ define isomorphisms since they satisfy 
\begin{align*}
E_{ev}Rv= v \hbox{ for all } v\in \cD_{per,ev}(X), \quad \hbox{and} \quad
RE_{ev}v= v \hbox{ for all } v\in \cD_{N}(X).
\end{align*}
If $X$ is reflexive, 
then for $p\in(1,\infty)$ and $s\in [0,\infty)$, we define the spaces
\begin{align*}
H^{s,p}_{N}(-h,0;X) := \overline{\cD_{N}}^{\norm{\cdot}_{H^{s,p}(-h,0;X)}}\quad \hbox{and} \quad
H^{-s,p}_{N}(-h,0;X) := H^{s,p'}_{N}(-h,0;X')'
\end{align*}
identifying $X$ with its bi-dual $X''$.
The completion has been taken with respect to the  $H^{s,p}(-h,0;X)$-norm of the $X$-valued Bessel potential spaces, compare e.g. \cite[Subsection VII.1.2]{Amann2019} for their definition. For $s \in \IN$ these are the classical Sobolev space, i.e.,   $H^{s,p}=W^{s,p}$. Using Bochner spaces, cf. e.g. \cite[Section 1.1]{ABHN}, one can define consistently $L^{p}(-h,0;X)$ even for $p\in [1,\infty]$, and one identifies
$L^p(-h,0;L^p(G;\cH))=L^p(\cD;\cH)$ for $p\in [1,\infty]$ using Fubini's theorem. 
Then, we define
\begin{eqnarray*}
H_z^{s,p}(H^{m,q}_{x,y}(\cH)):= H_N^{s,p}((-h,0);H_{per}^{m,q}(G;\cH))  \quad \hbox{for } m,s \in \IR \hbox{ and } p,q\in (1,\infty),
\end{eqnarray*}
where we use the short hand notation $H_z^{s,p}H^{m,q}_{x,y}$ if $\cH\in \{\IR,\IC\}$ or if there is no ambiguity.

The spaces $\cD_{per,ev}(X)$ and $\cD_{N}(X)$ are 
equipped with the topology induced by the semi-norms given by $\sup_{z}\norm{\partial^{\alpha}_zv(z)}_{X}$ for $\alpha\in \IN_0$.
Hence, one defines the spaces of distributions $\cD_{per,ev}'$ and $\cD_{N}'$ as the respective topological duals.
For these the pull backed maps  of $E_{ev}$ and $R$ again induce isomorphisms.
On these spaces of distributions we define the map
\begin{align*}
J_z^{s}:=(1-\partial_z^2)^{s/2} \quad \hbox{for} \quad s\in \IR.
\end{align*}
If $X$ is a UMD space, then $J_z^{s}$ defines a map on the distributions on the torus which restricts to a map on $\cD_{per,ev}(X)$. In fact, this is a discrete Fourier multiplier with symbol $(1+k_z^2)^{s/2}$ for $k_z\in (\pi/h)\IZ$   which induces an isomorphism
\begin{align*}
J_z^{s}\colon H^{m,p}_N(-h,0;X)  \rightarrow H^{m-s,p}_N(-h,0;X), \quad s,m\in \IR, \hbox{ and } p\in (1,\infty)
\end{align*}
suppressing the dependence on $X$ and $m,r,p$ in this notation, cf. e.g. \cite[Section 2.1]{Amann2019}. Also, one has isomorphisms
\begin{align*}
J_H^{m}\colon L^{q}(G;\cH)  \rightarrow H^{m,q}_{per}(G;\cH), \quad
J_H^{m}v = (1-\Delta_H)^{m/2}v,
\quad m\in \IR, \hbox{ and } q\in (1,\infty)
\end{align*}
Then
\begin{align*}
\norm{v}_{H_z^{s,p}H^{m,q}_{x,y}}= \norm{ \norm{J^s_zJ_H^{m} v(\cdot,z)}_{L^{q}(G;\cH)}}_{L^{p}(-h,0)}\quad \hbox{for } m,s \in \IR \hbox{ and } p,q\in (1,\infty).
\end{align*}

\subsection{Solenoidal spaces and the hydrostatic Helmholtz projection}
Similarly to the Navier-Stokes equations, an appropriate framework for the primitive equations are  {\em hydrostatically  solenoidal vector fields} satisfying
\begin{align*}
\div_H \overline{v} = 0 \quad \hbox{where}\quad \overline{v}=\frac{1}{h}\int_{-h}^0 v(\cdot,\cdot, \xi)d\xi.
\end{align*}
Hence, for $\cH$ a separable Hilbert space, one defines for $p\in (1,\infty)$
\begin{align*}
\lpso &:= \overline{\{v\in \cD_{N}(-h,0; C^{\infty}_{per}(\overline{G});\cH)^2 \colon \div_H \overline{v} = 0\}}^{\norm{\cdot}_{L^{p}(\cD;\cH)^2}},
\end{align*}
cf. \cite{HK16} for the scalar valued case with $\cH\in \{\IR,\IC\}$ for which we use the abbreviation $L^p_{\overline{\sigma}}(\cD)$. 

Moreover,  there exists a continuous projection, the \textit{hydrostatic Helmholtz projection}, 
\begin{align}\label{eq:IP}
\IP\colon L^p(\cD;\cH)^2 \rightarrow \lpso, \quad \IP v = \tilde{v} + \IP_2 \overline{v},
\end{align}
where
\begin{align*}
\tilde{v} = v-\overline{v}, \quad \overline{v}=\frac{1}{h}\int_{-h}^0 v(\cdot,\cdot, \xi)d\xi,
\end{align*}
and $\IP_2$ is the actual two-dimensional Helmholtz projection on $G$ defined for periodic functions. The scalar valued case has been discussed in \cite{HK16,GGHHK17}. The more general $\cH$-valued case can be drawn back to the scalar case by considering an orthonormal basis $(e_n)$ of $\cH$, and then  the corresponding  statements follow componentwise from the scalar valued case. Note that $\IP$ is given by a discrete Fourier multiplier, cf. e.g. \cite[Section 2.1]{Hu2020} for the $L^2$-case, which extends to all spaces $H_N^{s,p}(-h,0;H_{per}^{r,p}(G;\cH))$ for $s,r\in \IR$ and $p\in (1,\infty)$. In particular the averaging in \eqref{eq:IP} can be replaced by a projection onto the $0th$ Fourier mode in $z$-direction.

So, to include solenoidal functions, we define for $p\in (1,\infty)$ and $s,r\in \IR$
\begin{align*}
X^{s,m}_{\overline{\sigma}, p}(\cH):=\{v\in H_{N}^{s,p}(-h,0;H^{m,p}_{per}(G;\cH))\colon  \IP v =v\},
\end{align*} 
and we set 
\begin{align*}
X^{s}(\cH)_{\overline{\sigma}, p}:= X^{s,0}_{\overline{\sigma}, p}(\cH), \quad \hbox{where} \quad X^{0}_{\overline{\sigma}, p}(\cH)=\lpso. 
\end{align*}
If $\cH\in \{\IR,\IC\}$ and if there is no ambiguity we shorten the notation to $X^{s,m}_{\overline{\sigma}, p}=X^{s,m}_{\overline{\sigma}, p}(\cH)$.
Note that here due to the periodicity 
\begin{align*}
\IP H_{N}^{s,p}(-h,0;H^{r,p}_{per}(G;\cH))= X^{s,m}_{\overline{\sigma},p}(\cH). 
\end{align*}
\begin{remark}[Geometry of the spaces]\label{rem:type2}
The spaces $H_{N}^{s,p}(-h,0;H^{m,p}_{per}(G;\cH))$ and $X^{s,m}_{\overline{\sigma},p}(\cH)$ are UMD spaces and for $p\geq 2$ they are of type $2$. This holds since the spaces $X^{s,m}_{\overline{\sigma},p}(\cH)\subset H_{N}^{s,p}(-h,0;H^{m,p}_{per}(G;\cH))$ are closed subspaces, and
the latter are isomorphic to
$L^{p}(-h,0;L^{p}(G;\cH))$ which has the respective properties, cf. e.g. \cite{Hyt2017}.
\end{remark}

\begin{remark}[Extension of operators from scalar to $\cH$-valued spaces]\label{rem:extensionH}
Note that by  \cite[I.8.24]{Ste93} bounded operators on $L_{\overline{\sigma}}^r(\cD)$ admit an extension to $L_{\overline{\sigma}}^r(\cD;\cH)$ with  identical norm and analogously for the other function spaces discussed above. 
\end{remark}

\subsection{The hydrostatic Stokes operator}\label{subsec:hydrostaticStokes}
Next, observe that the operators
\begin{align*}
\Delta_{z}v=\partial_z^2v  \quad\hbox{and}\quad \Delta_H v= (\partial_x^2+\partial_y^2)v
\end{align*}
are given via discrete Fourier multipliers on periodic test functions, and therefore  these extend to operators in $H^{s,p}_N(-h,0;H^{m,p}_{per}(G;\cH))^2$ with
\begin{align*}
\Delta_{z,N}v &= \Delta_z v, \quad v \in H^{s+2,p}_N(-h,0;H^{m,p}_{per}(G;\cH))^2,  \quad\hbox{and}\quad \\
\Delta_{H,N}v &=\Delta_H v, \quad v\in H^{s,p}_N(-h,0;H^{m+2,p}_{per}(G;\cH))^2,
\end{align*}
respectively. Note that the operators $\Delta_{z,N}$, $\Delta_{H,N}$, and $\IP$ commute since these operators are compatible with the extension and restriction operators $E_{ev}$ and $R$, and hence 
they are similar to the corresponding operators in the periodic spaces where they are given in terms of symbols. 
Thereby, the analysis reduces to the periodic setting, and this 
has been used for instance also in many works by Cao, Li and Titi, compare e.g. \cite{LiTiti2016} and the references therein. In particular Fourier series methods are available. 
This allows us to define on $X^{s,m}_{\overline{\sigma},p}(\cH)$ for $s\in \IR$ and $p\in (1,\infty)$ the \emph{hydrostatic Stokes operator}
\begin{align}\label{eq:def_As}
A_{p,\cH}^{(s,m)}  := -\IP \Delta_{z,N} - \IP  \Delta_{H,N} 
\quad \hbox{with } D(A_{p,\cH}^{(s,m)}) =X^{s+2,m}_{\overline{\sigma},p}(\cH) \cap X^{s,m+2}_{\overline{\sigma},p}(\cH).
\end{align}
For simplicity we set $A_{p,\cH}^{(s)}:= A_{p,\cH}^{(s,0)}$.
This is consistent with the definition in \cite[Section 4]{HK16} for the scalar case, since
\begin{align}\label{eq:strong_Stokes}
A_{p,\cH} := A_{p,\cH}^{(0)} = -\IP \Delta v, \quad D(A_p) :=  \{v\in H_{per,H}^{2,p}(\cD;\cH)^2 \cap \lpso\colon \restr{\partial_z v}{\Gamma_N} = 0\} ,
\end{align}
where 
the part of the boundary where Neumann boundary conditions are imposed is abbreviated by
\begin{align*}
\Gamma_N = \Gamma_{u}\cup \Gamma_{b}.
\end{align*}
We drop the subscript and write $A_{p}^{(s,m)}, A_p^{(s)}$, and $A_p$, respectively, if $\cH\in \{\IR,\IC\}$ or if there is no ambiguity.

\begin{remark}[Properties of the hydrostatic Stokes operator]\label{rem:hydrStokes}
Note that the hydrostatic Helmholtz projection and the operators $\Delta_{z,N}$ and $\Delta_{H,N}$ are resolvent commuting
since one can reduce this to the periodic setting. 
Note that this does not hold any more if other than Neumann boundary conditions on top and bottom are imposed, cf. e.g. \cite{GGHHK17}.	
Moreover, 
here the Kalton-Weis theorem on the sum of commuting operators, see \cite{KalWei01}, is applicable, giving the closedness of $A_{p,\cH}^{(s,m)}$
and that
\begin{align*}
\IP \Delta_{z,N} + \IP  \Delta_{H,N} 
=  \Delta_{z,N}\IP +   \Delta_{H,N}\IP = (\Delta_{z,N} +   \Delta_{H,N})\IP.
\end{align*}
In addition $A_{p,\cH}^{(s,l)}$ and $A_{p,\cH}^{(0,0)}$ are similar.
\end{remark}

\subsection{Maximal $L^q_t$-regularity and interpolations spaces}\label{subsec:maxreg_trace}
One setting for the Cauchy problem 
\begin{align*}
(\partial_t + A)v = f, \quad v(0)=v_0,
\end{align*}
assuming that $X_1,X_0$ are Banach spaces such that $X_1$ is densely embedded into $X_0$ and $A\in \mathcal{L}(X_1,X_0)$,
are time-weighted vector valued $L^q$- and Sobolev spaces. For $q \in (1,\infty)$, $\mu\in (1/q,1]$, $T\in (0,\infty]$ these are defined by  
\begin{align*}
L^{q}_{\mu}(0,T;X_i)&:=\{v \in L^1_{loc}(0,T;X_i): t^{1-\mu}v\in L^{q}(0,T;X_i)\} \mbox{ and } \\
H^{1,q}_{\mu}(0,T;X_i)&:=\{v\in L^{q}_{\mu}(0,T;X_i) \cap H_{loc}^{1,1}(0,T;X_i)\colon \partial_t v \in L^q_{\mu}(0,T;X_i)\} \quad \hbox{for }i \in \{0,1\},
\end{align*}
cf. \cite[Section 3.2.4]{PruessSimonett}. The 
time weighted maximal regularity space  is then
\begin{align*}
L^{q}_{\mu}(0,T;X_1)\cap 
H^{1,q}_{\mu}(0,T;X_0) \subset 
L^{q}_{\mu}(0,T;X_0).
\end{align*}
The natural trace space  $X_{\mu-1/q,q}$ is determined  by means of real interpolation spaces $X_{\theta,q}:=(X_0, X_1)_{\theta,q}$ for $\theta \in (0,1)$, see e.g. \cite[Section 3.3.4]{PruessSimonett}.

Here, for the operators $A_{p,\cH}^{(s,m)}$ 
we denote for fixed $T\in (0,\infty]$ the ground space and the time weighted maximal regularity space  
\begin{align*}
\Ee^{q}_{0,\mu}(X^{s,m}_{\overline{\sigma},p}(\cH))&:= L^q_\mu(0,T; X^{s,m}_{\overline{\sigma},p}(\cH)) \\
\Ee^{q}_{1,\mu}(X^{s,m}_{\overline{\sigma},p}(\cH))&:=H_\mu^{1,q}(0,T; X^{s,m}_{\overline{\sigma},p}(\cH)) \cap L^q_\mu(0,T;D(A_{p,\cH}^{(s,m)}),
\end{align*}
which we abbreviate by $\Ee_0$ and $\Ee_1$, respectively, if there is no ambiguity. 
To formulate the following description of  trace spaces
we define the following Besov spaces.
For $s\geq 0$, $p,q\in(1,\infty)$, and a Banach space $X$ let
\begin{align*}
B^{s}_{pq,per}(G;\cH) &:= \overline{C^{\infty}_{per}(\overline{G};\cH)}^{\norm{\cdot}_{B^{s}_{pq}(G;\cH)}}, \\ 
B^{s}_{pq,per, ev}(-h,h;X)&:=\overline{\{u\in C^{\infty}[-h,h];X) \hbox{ periodic and even}  \}}^{\norm{\cdot}_{B^{s}_{pq,}(-h,h;X)}},  \\ B^{s}_{pq,N}(-h,0;X)&:=B^{s}_{pq,per, ev}(-h,h;X)\vert_{(-h,0)},
\end{align*}
and for $s<0$ the corresponding spaces are defined via duality, 
where $B^{s}_{pq}$ denotes Besov spaces which are defined as restrictions of Besov spaces on the whole space, cf.  \cite[Definitions 3.2.2]{Triebel}.

\begin{proposition}[Bounded $H^\infty$-calculus]\label{prop:Hinfty}
For $s,m\in \IR$ the operator $A_{p,\cH}^{(s,m)}+\nu$ admits for any $\nu >0$ a bounded $H^\infty$-calculus of angle zero. In particular, $A_{p,\cH}^{(s,m)}$ has the property of deterministic  maximal $L^q_{\mu}$-regularity on any finite time interval $(0,T)$ for $T\in (0,\infty)$ and for $p,q\in (1,\infty)$, $\mu\in (1/q,1]$, that is the solution operator for
\begin{align*}
(\partial_t+A_{p,\cH}^{(s,m)})v=f, \quad v(0)=0
\end{align*}
is an isomorphism 
\begin{align*}
\Ee^{q}_{0,\mu}(X^{s,m}_{\overline{\sigma},p}(\cH))\rightarrow
\{v\in \Ee^{q}_{1,\mu}(X^{s,m}_{\overline{\sigma},p}(\cH))\colon v(0)=0\}, \quad f \mapsto v.
\end{align*}
Moreover, for $\theta\in (0,1)$ and $q\in (1,\infty)$
\begin{align*}
D((A_{p,\cH}^{(s,m)})^\theta)
&=X^{s+2\theta,m}_{\overline{\sigma},p}(\cH) \cap X^{s,2\theta+m}_{\overline{\sigma},p}(\cH) \hbox{ and}\\
(X^{s,m}_{\overline{\sigma},p}(\cH), D(A_{p,\cH}^{(s,m)}))_{\theta, q}
&=
\IP B^{s+2\theta}_{pq,N} (-h,0; H^{m,p}(G;\cH)^2) \cap 
\IP H^{s,p}_{N} (-h,0; B^{2\theta+m}_{pq,per}(G;\cH)).
\end{align*}
\end{proposition}

\begin{proof}
It was shown in \cite[Theorem 3.1]{GGHHK17} that $A_p+\nu$ admits for any $\nu >0$ a bounded $H^\infty$-calculus of angle zero  in $L^p_{\overline{\sigma}}(\cD)$ and by Remark~\ref{rem:extensionH} this carries over to $\lpso$, and then by  similarity -- cf. Remark~\ref{rem:hydrStokes} --  to $A_{p,\cH}^{(s,m)}+\nu$ for $s,m\in \IR$.   

The domains of the fractional powers are hence given via  the complex interpolation functor $[\cdot,\cdot]_\theta$ as
	\begin{align*}
	D((A_{p,\cH}^{(s,m)})^\theta)&=[X^{s,m}_{\overline{\sigma},p}(\cH), D(A_{p,\cH}^{(s,m)})]_{\theta}\\
	&=[X^{s,m}_{\overline{\sigma},p}(\cH), X^{s+2,m}_{\overline{\sigma},p}(\cH) \cap 
	X^{s,2+m}_{\overline{\sigma},p}(\cH)]_{\theta}\\
	&=[X^{s,m}_{\overline{\sigma},p}(\cH), X^{s+2,m}_{\overline{\sigma},p}(\cH)]_{\theta} \cap 
	[X^{s,m}_{\overline{\sigma},p}(\cH), X^{s,2+m}_{\overline{\sigma},p}(\cH)]_{\theta}\\
	&=X^{s+2\theta,m}_{\overline{\sigma},p}(\cH) \cap X^{s,2\theta+m}_{\overline{\sigma},p}(\cH).
	\end{align*}
	Here to show the last equality, one combines  
	retract and co-retract arguments with interpolation results for Bessel potential spaces, compare e.g. \cite{HHK2016} for a similar argument. The pre-ultimate equality follows from 
	\cite[Lemma 9.5]{EscherPruessSimonett2003}
	or \cite[Theorem 3.1]{GGS1993} which are applicable here since in \eqref{eq:def_As} the operators $A_{p,\cH}^{(s,m)}$ have been defined as the sum of two resolvent commuting operators with -- up to a shift -- bounded $H^\infty$-calculus of angle zero.

Similarly, one determines that
the real interpolation spaces are given 
in terms of Besov spaces as
\begin{align*}
(X^{s,m}_{\overline{\sigma},p}(\cH), X^{s+2,m}_{\overline{\sigma},p}(\cH))_{\theta,q} &= \IP B^{s+2\theta}_{pq,N} (-h,0; H^{m,p}(G;\cH)^2), \quad \hbox{and}\\
(X^{s,m}_{\overline{\sigma},p}(\cH), X^{s,2+m}_{\overline{\sigma},p}(\cH))_{\theta,q} &=
\IP H^{s,p}_{N} (-h,0; B_{pq,per}^{2\theta+m}(G;\cH)^2). \qedhere
\end{align*}
\end{proof}

\subsection{The hydrostatic Neumann map}\label{subsec:Neumann}
Consider 
the stationary hydrostatic Stokes equation 
\begin{align}\label{eq:PE_stationary}
\left\{
\begin{array}{rll}
- \Delta V + \nabla_H P_s   &   = 0, &\text{ in } \cD,  \\
\mathrm{div}_H \overline{V} & = 0, &\text{ in } \cD,\\
\int_{\cD}V & = 0, &
\end{array}\right.
\end{align}
subject  to the boundary conditions  for  given  $g$
\begin{align}\label{eq:bc_l_and_u}
\begin{split}
V, P_s  \hbox{ are periodic } \hbox{on } \Gamma_l, \quad \partial_z V = 0 \hbox{ on } \Gamma_b \quad \hbox{and} \quad \partial_z V = g \hbox{ on } \Gamma_{u}.
\end{split}
\end{align}
The boundary data will be taken from the $\cH$-valued Sobolev-Slobodeckij  spaces
\begin{align*}
W^{s,r}_{per}(\Gamma_{u};\cH):= \overline{C^{\infty}_{per}(\overline{\Gamma_{u}};\cH)}^{\norm{\cdot}_{W^{s,r}(\Gamma_{u};\cH)}}, 
\end{align*}
where $W^{s,r}(\Gamma_{u};\cH)$ is defined as restriction of $W^{s,r}(\IR^2;\cH)$. One has also $W^{s,r}_{per}(\Gamma_{u};\cH)=	B^{s}_{rr,per}(\Gamma_{u};\cH)$.

\begin{proposition}[Hydrostatic Neumann map]\label{prop:NeumannMap}
Let $\cH$ be a separable Hilbert space, $r\geq 2$, and $m\geq 0$.
For $g\in W^{1-1/r+m,r}_{per}(\Gamma_{u};\cH)^2$
there exist unique 
\begin{align*}
V&\in L^r_{\overline{\sigma}}(\cD;\cH)\cap H^{2+m,r}(\Omega;\cH)^2 \hbox{ with } \int_{\cD} V =0 \hbox{ and } \\
P_s&\in  H^{1,r}(G;\cH) \hbox{ with } \int_{G} P_s =0
\end{align*}
solving  
\eqref{eq:PE_stationary}
and \eqref{eq:bc_l_and_u}.
The Neumann map  given by
\begin{align*}
\mathcal{N}&\colon	W^{1-1/r+m,r}_{per}(\Gamma_{u};\cH)^2   \rightarrow L^r_{\overline{\sigma}}(\cD;\cH)\cap H^{2+m,r}(\Omega;\cH)^2, \quad g \mapsto V
\end{align*}
is continuous. In particular the following embedding is continuous
\begin{align*}
\mathcal{N}(	W^{1-1/r+m,r}_{per}(\Gamma_{u};\cH)^2) \hookrightarrow D(A_{r,\cH}^{(s,m)}) \quad \hbox{for } s \in [-3/2, 1/r-1).
\end{align*}	
\end{proposition}
For the Neumann map in the setting of diffusion equations, we refer to \cite{Amann1993}. 
For an abstract version of the Neumann map (there called Dirichlet operator)
we refer to \cite{BE20, CENN03, Gre87} and the references therein.

\begin{proof}[Proof of Proposition \ref{prop:NeumannMap}]
To rewrite the problem \eqref{eq:PE_stationary} and  \eqref{eq:bc_l_and_u},
we introduce the \emph{maximal hydrostatic Stokes operator} $A_{\max}$ 
in $L^r_{\overline{\sigma}}(\cD,\cH)$ for average free functions given by
\begin{equation*}
A_{\max} V := - \IP \Delta V , \quad 
D(A_{\max}) := \{ V \in W^{2,r}_{per}(\cD;\cH)^2 \cap L^r_{\overline{\sigma}}(\cD;\cH) \colon \partial_z V|_{\Gamma_{b}} = 0,   \int_{\cD} V =0  \} .
\end{equation*}
Further, consider the outer normal derivative on the upper boundary
\begin{align*}
L \colon D(A_{\max}) \to W^{1-1/r,r}_{per}(\Gamma_{u};\cH)^2, \quad f \mapsto (\partial_z f)(\cdot,\cdot,0).
\end{align*}
We point out that $A_{\max}|_{\ker(L)} = A_r|_{\{V\colon \int_{\cD}V=0\}}$ is the hydrostatic Stokes operator with homogeneous Neumann boundary conditions restricted to the average free functions.
Then, \eqref{eq:PE_stationary} with \eqref{eq:bc_l_and_u} is equivalent to 
\begin{equation}
\begin{cases}
A_{\max} V &= 0, \\
L V &= g .
\end{cases}
\label{eq:PE_stationary2}
\end{equation}

\begin{lemma}\label{lemma:L surjective} For $r\in (1,\infty)$
	the operator $L$
	is surjective. 
\end{lemma}
\begin{proof} 
	By \cite{Triebel} it follows using periodic extensions and restrictions that there is for $g \in W^{1-1/r,r}_{per}(\Gamma_u;\cH)^2$ 
	\begin{align*}
	G_g \in W^{2,r}_{per}(\cD;\cH)^2 \quad \hbox{such that}\quad \partial_z G_g(\cdot,\cdot,0) = g \hbox{ and } \partial_z G_g(\cdot,\cdot,-h) = 0.
	\end{align*}
	We set $\varphi := \div_H \overline{G} \in W^{1,r}(G;\cH)^2$. Now the Bogovskii operator $B$ on $G$ for periodic boundary conditions  yields a function $H_{\varphi} := B\varphi \in W^{2,r}(G;\cH)^2$ which we extend constantly in $z$-direction and with a slight abuse of notation still denote by $H_\varphi$. Then $H_\varphi \in W^{2,r}(\cD;\cH)^2$, and it satisfies $\partial_z H_\varphi = 0$, in particular $\partial_z H_\varphi(\cdot,\cdot,0) =0$ and $\partial_z H_\varphi(\cdot,\cdot,-h) = 0$, and $\overline{H_\varphi} = H_\varphi$, and therefore $\div_H(\overline{H_\varphi}) = \div_H(H_\varphi) = \varphi$. Hence $F := G_g - H_\varphi \in W^{2,r}(\cD;\cH)^2$ satisfies 
	\begin{equation*}
	\div_H(\overline{F}) = \div_H(\overline{G_g}) - \div_H(\overline{H_\varphi})
	= \varphi - \varphi = 0
	\end{equation*}
	and $\partial_z F(\cdot,\cdot,-h) = 0$, i.e. $\tilde{F}:=F - \int_{\cD}F\in D(A_{\max})$. 
	Finally, we obtain
	\begin{equation*}
	L \tilde{F} =  \partial_zF(\cdot,\cdot,0) = \partial_z G_g(\cdot,\cdot,0)  - \partial_z H_\varphi(\cdot,\cdot,0) 
	= g. \qedhere
	\end{equation*}
\end{proof}
The operator $A_r|_{\{V\colon \int_{\cD}V=0\}}= A_{\max}|_{\ker(L)}$ is boundedly invertible, which follows from \cite{GGHHK17}. 
Further, $A_{\max}$ is closed and $L$ is relatively $A_{\max}$-bounded by the trace theorem. Hence
\begin{align*}
\binom{A_{\max}}{L} \colon D(A_{\max}) \to \{V\in L^r_{\bar{\sigma}}(\cD;\cH)\colon \int_{\cD} V=0\} \times W_{per}^{1-1/r,r}(\Gamma_{u};\cH)^2
\end{align*}
is closed. Finally, by Lemma~\ref{lemma:L surjective} the operator $L \colon D(A_{\max}) \to W_{per}^{1-1/r,r}(\Gamma_{u};\cH)^2$ is surjective. 
Now \cite[Lemma~2.2]{CENN03} implies that for every $g \in W_{per}^{1-1/r,r}(\Gamma_{u};\cH)^2$ there exists a unique solution $V$ of \eqref{eq:PE_stationary2} and that the Neumann map
\begin{equation*}
\mathcal{N} \colon W^{1-1/r,r}_{per}(\Gamma_{u};\cH)^2   \rightarrow
\{V\in L^r_{\bar{\sigma}}(\cD;\cH)\colon \int_{\cD} V=0\}, \quad g \mapsto V 
\end{equation*}
is continuous.
For $m>0$, note that the solution $V$ of \eqref{eq:PE_stationary2} satisfies in particular $\Delta V = \nabla_H P_s$ and
\begin{align*}
\partial_z V|_{\Gamma_{u}} = g, \quad 
\partial_z V|_{\Gamma_{b}} = 0, \quad 
\int_{\cD} V = 0, \quad \hbox{and} \quad
V \text{ periodic } \text{on } \Gamma_l.
\end{align*}
Using the boundary  and the divergence free conditions it follows 
by taking the vertical average and applying $\div_H$, 	
compare also  \cite{GGHHK20},  that
\begin{align*}
-\div_H\Delta_H\overline{V}-\tfrac{1}{h}\div_Hg=-\tfrac{1}{h}\div_Hg=\Delta_H P_s,
\end{align*}
and hence
the pressure term $\nabla_H P_s$ is given by
\begin{equation*}
\nabla_H P_s = \frac{1}{h} \nabla_H \Delta_H^{-1} \div_H g .
\end{equation*}
Now $g \in W_{per}^{1-1/r+m,r}(\Gamma_u,\cH)^2$ implies $\nabla_H P_s \in W_{per}^{1-1/r+m,r}(\Gamma_u,\cH)^2$. Therefore the second claim follows by regularity theory of the Poisson equations.

For the last claim 
one uses that 
\begin{align*}
H^{s+2,r}(-h,0;X)= H_{z,N}^{s+2,r}(-h,0;X)\quad \hbox{for} \quad 0\leq 2+s<1+1/r.
\end{align*}
To include negative differentiability consider  for $1/r+1/r'=1$ 
\begin{align*}
H^{|s+2|,r'}(-h,0;X)= H_{z,N}^{|s+2|,r'}(-h,0;X)\quad \hbox{for} \quad 0\leq |2+s|<1+1/r'=2-1/r,
\end{align*}
and then by duality 
\begin{align*}
H^{s+2,r}(-h,0;X)= H_{z,N}^{s+2,r}(-h,0;X)\quad \hbox{for} \quad -2+1/r<2+s< 0.
\end{align*}
Moreover, 
\begin{align*} 
H^{2+m,r}(\cD;\cH)^2  = H^{2+m,r}(-h,0; L^r(G;\cH))\cap  L^{r}(-h,0; H^{2+m,r}(G;\cH)),
\end{align*}
and by the mixed derivative theorem  \cite[Corollary 4.5.10]{PruessSimonett} applied to $(-\Delta_H)^{(2+m)/2}$ and $(-\partial_z^2)^{(2+m)/2}$ with $\theta=2/(2+m)$ 
\begin{align*}
H^{2+m,r}_{z,N}(-h,0; L^r(G;\cH))\cap  L^{r}(-h,0; H^{2+m,r}_{per}(G;\cH)) \hookrightarrow 
H^{2}_{z,N}(-h,0; H^{m,r}_{per}(G;\cH)).
\end{align*}
Hence, since $2+s<1+1/r<2$, that is $s<1/r-1<0$,  and $L^r_z\hookrightarrow H^{s,r}_z$
\begin{multline*}
H^{2+m,r}(\cD;\cH)^2 \cap X_{\overline{\sigma},p}^{s,m}(\cH) 
\hookrightarrow \\
H_{z,N}^{2+s,r}(-h,0; H^{m,r}(G;\cH))\cap H_{z,N}^{s,r}(-h,0; H_{per}^{2+m,r}(G;\cH))\cap X_{\overline{\sigma},p}^{s,m}(\cH)
= D(A_{r,\cH}^{(s,m)}).
\qedhere
\end{multline*}  
\end{proof}

\section{Linear stochastic theory}\label{sec:prel_stoch}

\subsection{Cylindrical Wiener processes and stochastic convolutions}\label{subsec:Wiener}

Let $(\Omega, \mathcal{A},P)$ be a probability space with a filtration $\mathcal{F}=(\mathcal{F}_t)_t$. An 
$\mathcal{F}$-cylindrical  Wiener process (or cylindrical Brownian motion) on a separable Hilbert space $\cH$ is a 
bounded linear operator 
$$W_{\cH}\colon L^2((0,\infty);\cH)
\rightarrow L^2(\Omega)$$ 
such that for all $f,g\in \cH$ and $0\leq t\leq t'$:
\begin{itemize}
\item[(a)] The random variable $W_{\cH}(t)f:=W_{\cH}(\mathds{1}_{[0,t]}\otimes f)$ is centered Gaussian and $\mathcal{F}_t$-measurable;
\item[(b)] $\E[W_{\cH}(t')f\cdot W_{\cH}(t)g]= t \left< f,g \right>_\cH$;
\item[(c)] The random variable $W_{\cH}(t')f-W_{\cH}(t)f$ is independent of $\mathcal{F}_t$.
\end{itemize}
From now on we fix the separable Hilbert space $\cH$ and the filtration $\cF$.
Let $(e_n)_n$ be an orthonormal 
basis   of $\cH$, then $\beta_n(t):=W_{\cH}(t)e_n$ is a standard $\mathcal{F}$-Brownian motion, and we 
have the representation 
\begin{align*}
W_{\cH}(t)\colon \cH \to L^2(\Omega), \quad
W_{\cH}(t)f = \sum_{n=1}^\infty \beta_n(t) \left< f,e_n \right>_\cH,
\end{align*}
and $(W_{\cH}(t))_{t\geq 0}$ 
defines a family of linear operators called the cylindrical  Wiener process $W_{\cH}$. 

For the definition of the stochastic integral with respect to $W_{\cH}$ we refer to the paper by van Neerven, Veraar and Weis \cite{NeervenVeraarWeis}.
Mostly, the stochastic integral is defined for functions taking values in the space of 
$\gamma$-radonifying operators $\gamma(\cH;X^{s,m}_{\overline{\sigma},p})$,  see \cite[Chapter 9]{Hyt2017} for more details. Here, we will use instead $X^{s,m}_{\overline{\sigma},p}(\cH)$ similar to 
the $L^p$-setting in \cite{NeervenVeraarWeis}, because the spaces $\gamma(\cH;X^{s,m}_{\overline{\sigma},p})$ and $X^{s,m}_{\overline{\sigma},p}(\cH)$
are isomorphic by \cite[Theorem 9.3.6]{Hyt2017}, 
\begin{align}\label{eq:gamma}
\gamma(\cH;X^{s,m}_{\overline{\sigma},p}) \cong  X^{s,m}_{\overline{\sigma},p}(\cH) \quad \hbox{for } p\in (1,\infty), \hbox{ and } s,m\in \IR.
\end{align}
We will use the latter spaces throughout the paper where both ways to define the stochastic integral are equivalent.



Consider the  hydrostatic Stokes equations in $X^{s,m}_{\overline{\sigma},r}$ with $\cF$-adapted stochastic forcing $H_f$ with values in $X^{s,m}_{\overline{\sigma},r}(\cH)$
\begin{align}\label{eq:stoch}
\d Z_f(t) + A_r^{(s,m)} Z_f(t) \dt = H_f(t) \dW_\cH(t), \quad Z(0)=Z_0.
\end{align}
A \emph{strong solution} is given 
by the stochastic convolution
\begin{align}\label{eq:stochconvolution}
Z_f(t):=e^{-tA_r^{(s,m)}}Z_0 + \int_0^t e^{-(t-s)A_r^{(s,m)}} H_f(s) \dW(s).
\end{align}
if certain regularity conditions made precise in Proposition 
\ref{theorem:stochlinex} below hold.
By Proposition~\ref{prop:Hinfty} and Remark~\ref{rem:type2}, we can apply here the theory of stochastic maximal 
$L^r$-regularity developed by van Neerven, Veraar and Weis, cf.  \cite{NeervenVeraarWeis, NeervenVeraarWeis2} to define and estimate $Z_f$ 
given in
\eqref{eq:stochconvolution}.
For the case discussed here, we  summarize the result as follows applying the results from \cite[Section 3]{AgrestiVeraar}. 
Note that the time weights  in there translate to the ones used here (which follows the notation in \cite{PruessSimonett}) by 
$\kappa/q =1-\mu$, where $\kappa \in [0, q/2-1)$ that is $\mu = 1-\kappa/q$ and $\mu \in (1/2+1/q,1]$.

\begin{proposition}[Stochastic maximal regularity for the hydrostatic Stokes equations]\label{theorem:stochlinex} 
Let $0<T<\infty$, $s,m\in \IR$, $r,q\geq 2$ with $r>2$ if $q\neq 2$, and $\mu \in (1/2+1/q,1]$ with $\mu=1$ if $p=q=2$. Then for any strongly $\mathcal{F}_0$-measurable
\begin{align*}
Z_0\colon\Omega\to (X^{s,m}_{r, \overline{\sigma}},D((A_r^{(s,m)}))_{\mu-1/q,q} \quad \hbox{and}\quad H_f\in  L^q(\Omega;L^q_{\mu}(0,T;D((A_{r,\cH}^{(s,m)})^{1/2}) )
\end{align*} being $\mathcal{F}$-adapted,  
the stochastic convolution \eqref{eq:stochconvolution} is well-defined, $Z_f$ given by \eqref{eq:stochconvolution} is $\mathcal{F}$-adapted and defines the  unique solution to \eqref{eq:stoch} satisfying pathwise
\begin{align*}
Z_f &\in H_{\mu}^{\theta,q}\big( 0,T; D((A_r^{(s,m)})^{1-\theta})\big)
\cap C\big( [0,T]; (X^{s,m}_{r, \overline{\sigma}},D(A_r^{(s,m)}))_{\mu-1/q,q}\big) \quad \hbox{for any } \theta\in[0,1/2).
\end{align*}
\end{proposition}
Note that here $H_f$ is defined in the general $\cH$-valued setting, where $\cH$ is a separable Hilbert space and we have used the identification \eqref{eq:gamma}, while all other functions have values in finite dimensional spaces.


\begin{remark}[$L^q$-estimates in the probability space]\label{rem:thm_stoch}
Assuming in Proposition~\ref{theorem:stochlinex} additionally $Z_0\in L^q(\Omega;  (X^{s,m}_{r, \overline{\sigma}},D((A_r^{(s,m)}))_{\mu-1/q,q} )$, it follows that in Proposition~\ref{theorem:stochlinex}  even 
\begin{align*}
Z_f \in L^q\big(\Omega;H_{\mu}^{\theta,q}\big( 0,T; D(A_r^{1/2-\theta})\big) \cap C\big( [0,T]; (X^{s,m}_{r, \overline{\sigma}},D((A_r^{(s,m)}))_{\mu-1/q,q} \big)\big), \quad \theta\in[0,1/2),
\end{align*}
and the corresponding stochastic  maximal $L^q$-regularity estimate holds, see \cite[Proposition 3.10]{AgrestiVeraar}. 
However, the pathwise regularity result of Proposition~\ref{theorem:stochlinex} is sufficient for our purpose to construct pathwise solution.  
\end{remark}

\subsection{Rewriting the  stochastic boundary condition  as a forcing term}\label{subsec:boundary}
To give a precise definition of a solution to \eqref{eq:pe_stoch} with boundary conditions \eqref{eq:bc_l} - \eqref{eq:bcstoch}, we assume that there are functions
\begin{align}\label{eq:hbg}
g \colon \Gamma_{u}\rightarrow \cH^2\quad \hbox{and}\quad  h_b\colon  \Gamma_u\times (0,T) \rightarrow \IR \quad \hbox{with } \int_{\Gamma_u}h_b(\cdot)g=0,
\end{align}
and a cylindrical Wiener process $W_b$ defined on $\cH$ with respect to the
filtration $\mathcal{F}$ such that the noise term can be written given an orthonormal basis $(e_n)$ of $\cH$ as
\begin{align*}
\omega(t) = \sum_{n=1}^{\infty} <g,e_n>  W_b(t) e_n=W_b(t)g.
\end{align*}
Examples of such  $\omega$ have been discussed also in \cite[Example 4.6]{SchnaubeltVeraar2011}. 
For $r \geq 2$ and $s,m\in \IR$ let
\begin{align*}
h_b(t)g \in W^{1-1/r+m,r}_{per}(\Gamma_{u};\cH)^2\quad \hbox{for }  t\in (0,T). 
\end{align*}
Considering the equation
\begin{align*}
\d Z_b(t) + A_r^{(s,m)} Z_b(t) \dt =0, \quad Z_b(0)=0,
\end{align*}
subject to \eqref{eq:bc_l}, \eqref{eq:bc_Neumann}, and the inhomogeneous stochastic boundary conditions
\begin{align*}
\partial_z Z_b = h_b \partial_t \omega \hbox{ on } \Gamma_{u} \times (0,T),
\end{align*}
then
we call $Z_b$ a \emph{strong pathwise 
solution}, if
\begin{align}\label{eq:Zb}
Z_b(t):= A_r^{(s,m)} \int_0^t e^{-(t-s)A_r^{(s,m)}} [\mathcal{N} h_b(t)g] \dW_b(t),
\end{align}
where $\mathcal{N}$ denotes the hydrostatic Neumann operator defined in Subsection \ref{subsec:Neumann}
and the stochastic integral in defined as the one in \eqref{eq:stochconvolution}, compare also \cite[Section 13]{DZ96} for an analogous definition. Here we identify similarly to \eqref{eq:gamma}
\begin{align*}
\gamma(\cH;W^{s,r}_{per}(\Gamma_u)) \cong  W^{s,r}_{per}(\Gamma_u;\cH) \quad \hbox{for } r\in (1,\infty), \hbox{ and } s\in \IR.
\end{align*}

\begin{proposition}[Pathwise regularity of the boundary forcing]\label{prop:reg_Zb}
Let $0<T<\infty$, $s,m\in \IR$, $r,q\geq 2$ with $r>2$ if $q\neq 2$, and $\mu \in (1/2+1/q,1]$ with $\mu=1$ if $p=q=2$, and $$h_b(t)g \in L^q_{\mu}(0,T;W^{1-1/r+m,r}(\Gamma_u;\cH)^2),$$ then  pathwise
\begin{align*}
Z_b &\in L^q_{\mu}(0,T; D((A^{(s,m)}_r)^{1/2}) \quad \hbox{for } m\geq 0 \hbox{ and }s\in[-3/2, 1/r-1)
\end{align*}
\end{proposition}

\begin{proof}
By Proposition~\ref{prop:NeumannMap}, where the restriction on $s$ enters, one has
\begin{align*}
[\mathcal{N} h_b(t)g] \in L_{\mu}^q(0,T;H_{per,H}^{2,r}(\cD;\cH)\cap L^{r}_{\overline \sigma}(\cD;\cH))\subset
L_{\mu}^q(0,T;D(A_{r,\cH}^{(s,m)})).
\end{align*}	
Hence,
\begin{equation*}
A_r^{(s,m)}\int_0^t e^{-(t-s)A_r} [\mathcal{N} h_b(t)g] \dW_b(t) = \int_0^t  e^{-(t-s)A_r} A_{r,\cH}^{(s,m)}[\mathcal{N} h_b(t)g] \dW_b(t).
\end{equation*}
By the stochastic maximal regularity -- similarly to the situation in
Proposition~\ref{theorem:stochlinex} -- the spatial regularity increases by an order of $1/2$, and the claim follows. 
\end{proof}

\section{Local pathwise well-posedness}\label{sec:main}

\subsection{Notion of solution and main result}

Now we are in the position to adapt 
the approach by Da Prato and Zabczyk, compare~\cite[Section 13]{DZ96}, to define a notion of solution to the stochastic primitive equations with \emph{inhomogeneous} stochastic boundary conditions.
To this end, suppose that $V$ solves 
the primitive equations with inhomogeneous stochastic boundary conditions
\eqref{eq:pe_stoch}, and $Z_b$  be given by \eqref{eq:Zb}, that is,
\begin{align*}
Z_b(t):= A_r^{(s,m)} \int_0^t e^{-(t-s)A_r^{(s,m)}} [\mathcal{N} h_b(t)g] \dW_b(t).
\end{align*}
Then we set
\begin{align}
V_b:=V-Z_b. 
\end{align}
Applying the hydrostatic Helmholtz projection $\IP$, this solves the stochastic differential equation
\begin{align}\label{eq:sol_strong2}
&\d V_b+ A V_b \dt =F(V_b +Z) \dt + H_f \d W, \quad V_b(0)=V_0.
\end{align}
%
which is  subject to  the now \emph{homogeneous} boundary conditions  
\begin{align*}
V_b \hbox{ are periodic } \hbox{on } \Gamma_l \times (0,T),  
\quad \hbox{and} \quad \partial_z V_b = 0 \hbox{ on } (\Gamma_{u} \cup \Gamma_{b}) \times (0,T).
\end{align*}
Here we defined the bilinear map $F(\cdot,\cdot)$ by 
\begin{align}\label{eq:F}
F(v,v') := \IP(v \cdot \nabla_H v' + w(v)\partial_z v').
\end{align} 
As discussed in detail by Da Prato and Zabczyk in~\cite[Section 13]{DZ96}, if the solution $V$  were sufficiently regular, then it would constitute a solution to the actual the inhomogeneous problem.  This motivates the following notion of solution:
We call $V=V_b+Z_b$  a \emph{pathwise $z$-weak solution} to  \eqref{eq:pe_stoch} on $(0,T)$ for $T>0$ subject to  the inhomogeneous stochastic boundary condition  \eqref{eq:bc_Neumann}--\eqref{eq:bcstoch} 
with initial condition $V_0$ 
if $Z_b$ is given by \eqref{eq:Zb},  $V_0\colon \Omega \rightarrow X^{-1}_{\overline{\sigma},p}$,  
\begin{align*}
V_b \colon \Omega \rightarrow L^{q}_{\mu}\big( 0,T; D(A^{(-1)}_p)\big)\cap C([0,T];X^{-1}_{\overline{\sigma},p}),
\end{align*}
$V_b$ is adapted, and these solves pathwise the equation \eqref{eq:sol_strong2}. Note that then the surface pressure $P_s$ can be recovered from $V_b$.

Given $Z_f$ as in  \eqref{eq:stochconvolution}
which solves 
\begin{align*}
\d Z_f(t) + A Z_f(t) \dt = H_f(t) \dW(t),  \quad Z(0)=Z_0,
\end{align*}
one considers the difference
\begin{align*}
v:= V_b-Z_f = V-Z\quad \hbox{with} \quad Z:= Z_f+Z_b,
\end{align*}
where the stochastic term cancels out and  which therefore solves the \emph{deterministic} equation 
\begin{align}\label{eq:PE_stoch_z}
\partial_t v + A v =  F(v+Z,v+Z), \quad v(0)=v_0, \quad \hbox{where } v_0=V_0-Z_0.
\end{align}
This will be solved in the following framework.
Recall from Section~\ref{sec:prel} that for $p\in (1,\infty)$ and $s,m\in \IR$
\begin{align*}
X^{s,m}_{\overline{\sigma}, p}(\cH)&=\{v\in H_{N}^{s,p}(-h,0;H^{m,p}_{per}(G;\cH))\colon  \IP v =v\},\\
X^{-1}_{\overline{\sigma}, p}(\cH)&=\{v\in H_{N}^{-1,p}(-h,0;L^{p}_{per}(G;\cH))\colon  \IP v =v\},\\
A_{p,\cH}^{(s,m)} &= \IP \Delta_{z,N} + \IP  \Delta_{H,N} 
\quad \hbox{with } D(A_{p,\cH}^{(s,m)}) =X^{s+2,m}_{\overline{\sigma},p}(\cH) \cap X^{s,m+2}_{\overline{\sigma},p}(\cH),\\
A_{p,\cH}^{(-1)} &= \IP \Delta_{z,N} + \IP  \Delta_{H,N} 
\quad \hbox{with } D(A_{p,\cH}^{(-1)}) =X^{1,0}_{\overline{\sigma},p}(\cH) \cap X^{-1,2}_{\overline{\sigma},p}(\cH),
\end{align*} 
where we omit $\cH$ in the notation when $\cH\in \{\IC,\IR\}$ and we identify the $\cH$-valued spaces with the $\gamma$-radonyfying operators, cf. \eqref{eq:gamma}.  
We are now in the position to formulate the main result of this article.

\begin{theorem}[Local $z$-weak pathwise well-posedness]\label{thm:main1} 
Let $0<T<\infty$, $p\in (3,4]$,  $\varepsilon>0$ sufficiently small,  $\mu:=5/8+3/4p+3\varepsilon/8$ and $q>2$ such that $\mu\in (1/2+1/q,1]$ and $\sigma=(1+\mu)/2 \in (1/2+1/2q,1]$.
Let be given a
strongly $\mathcal{F}_0$-measurable
\begin{align*}
Z_0\colon\Omega\to (X^{-1}_{p, \overline{\sigma}},D(A_p^{(-1)}))_{\sigma-1/2q,2q} \quad \hbox{and}\quad H_f\in  L^{2q}(\Omega;L^{2q}_{\sigma}(0,T;D((A_{p,\cH}^{(-1)})^{1/2}) )
\end{align*} being $\mathcal{F}$-adapted,
\begin{align*}
h_b(t)g \in L^{2q}_{\sigma}(0,T;W^{1+1/p+\varepsilon,p}(\Gamma_u;\cH)^2), \quad \hbox{and} \quad v_0\in (X^{-1}_{\overline{\sigma}, p},D(A_p^{(-1)}))_{\mu-1/q,q}.
\end{align*}
Then there exists a unique, local in time $z$-weak pathwise solution to the stochastic primitive equations~\eqref{eq:pe_stoch} subject to \eqref{eq:bc_l}--\eqref{eq:bcstoch} 
\begin{align*}
V=v+Z_f +Z_b. 
\end{align*}
Here,
$Z_f$ is given by \eqref{eq:stochconvolution} and $Z_b$ by \eqref{eq:Zb} and  pathwise for   $s=1/2-1/p-\varepsilon/2$, and $\theta\in[0,1/2)$,
\begin{align*}
Z_f&\in H_{\sigma}^{\theta,2q}\big( 0,T; D((A_p^{-1})^{1-\theta})\big) \cap C\big( [0,T]; (X^{-1}_{p, \overline{\sigma}},D(A_p^{(-1)}))_{\sigma-1/2q,2q} \big),   \\
Z_b &\in
H_{\sigma}^{\theta,q}\big( 0,T; D((A_p^{(s,2/p+\varepsilon)})^{1/2-\theta})\big) \cap C\big( [0,T]; (X^{s,m}_{r, \overline{\sigma}},D(A_p^{(s,2/p+\varepsilon)}))_{\sigma-1/2q,2q} \big),
\end{align*}
and there is a $T^\ast\in (0, T]$ such that $v$ is a deterministic function with	
\begin{align*}
v&\in H^{1,q}_{\mu}(0,T^\ast;X_{\overline{\sigma},p}^{-1}) \cap L^{q}_{\mu}(0,T^\ast;D(A_p^{(-1)})),  
\end{align*} 
which 
depends continuously on $v_0$.
\end{theorem}

\begin{remark}
All terms in the equation~\eqref{eq:Vb}, reformulated  in  \eqref{eq:PE_stoch_z} below, can be made sense of within these regularity classes as the estimates in Lemma~\ref{lemma:nonlin_I} and Lemma~\ref{lemma:nonlin_II} show for the non-linear term, and the linear term is given via the duality pairing
\begin{align*}
\langle \partial_z v, \partial_z \phi\rangle + \langle \Delta_H v, \phi\rangle, \quad \phi \in \IP H_z^{1,p'}L^{p'}_{xy}\cap \IP H_z^{-1,p'}H^{2,p'}_{xy}.
\end{align*}  
The restriction $p\in (3,4]$ is a technical assumption due to  
the non-linear estimates, and the regularity of $Z_f$ and $Z_b$. On the one hand  certain embeddings 
which require $p>3$ are employed, cf. \eqref{eq:pge3} below. On the other hand, one has to assure that the embedding of the range of the Neumann map into a certain operator domain, see Proposition~\ref{prop:NeumannMap}, remains valid, where the larger $p$ the less derivatives are available, cf. \eqref{eq:pleq4} below. Balancing  these two requirements leads us here  to the range $p\in (3,4]$. 
\end{remark}

\begin{remark}[Blow up criteria]
For the deterministic part $v$ of the pathwise solution solving \eqref{eq:PE_stoch_z} below, there are blow-up criteria of Serrin type in maximal regularity spaces: Let $T^{\max}$ be the maximal existence time of $v$, then by \cite[Theorem 2.4]{PruessSimonettWilke}
	\begin{enumerate}
		\item $v\in L^q(0,T';D((A^{-1})^{\mu}))$ for all $T^{'}<T^{\max}$;
		\item if $T^{\max}<\infty$, then  $v\notin L^q(0,T^{\max};D((A_p^{(-1)})^{\mu}))$,
	\end{enumerate}  
	where $\mu$ is the time-weight from Theorem~\ref{thm:main1}.
Blow-up criteria in the stochastic -- not only pathwise setting -- are discussed in \cite[Section 4]{AgrestiVeraar}.
\end{remark}

\subsection{Non-linear estimates}\label{sec:proof}
Recall that, similarly to \cite[Section 5]{HK16}, we defined  
in \eqref{eq:F} the bilinear map $F(\cdot,\cdot)$ by 
\begin{align*}
F(v,v') := \IP(v \cdot \nabla_H v' + w(v)\partial_z v'),
\end{align*}
and set $F(v):=F(v,v)$. 
Using the product rule and that $(v,w(v))$ is divergence free one rewrites 
$F(v,v')
= \IP \partial_z(w v') + \IP \div_H(v\otimes v'),$
and hence 
\begin{align*}
F(v,v')=F_z(v,v')+ F_H(v,v')
\end{align*} 
with
\begin{align*} 
F_z(v,v'):=\IP(\partial_z (w(v) v')) 
\quad\hbox{and}\quad 
F_H(v,v')=\IP \div_H(v\otimes v'). 
\end{align*}
Moreover, we set for brevity
\begin{align*}
X_{\beta}:=D((A_p^{(-1)})^{\beta})
\quad \hbox{for } \beta \geq 0.
\end{align*}

\begin{lemma}[Estimate on the non-linearity part I]\label{lemma:nonlin_I}
Let $p\in (1,\infty)$, then there exists a constant $C>0$ depending only on $h,p$ such that for $\beta_{z}=1/2+3/4p$ 
\begin{align*}
\norm{F_z(v,v')}_{X^{-1}_{\overline{\sigma},p}} &\leq C 
\norm{v}_{X_{\beta_{z}}}\norm{v'}_{X_{\beta_{z}}}, \\
\norm{F_z(v,Z')}_{X^{-1}_{\overline{\sigma},p}}
&\leq C
\norm{v}_{X_{\beta_{z}}} \norm{Z'}_{H^{1/2p,p}_z H^{1/p,p}_{xy}}, \\
\norm{F_z(Z,v')}_{X^{-1}_{\overline{\sigma},p}}
&\leq C \norm{Z}_{H^{-1+1/2p,p}_zH^{1+1/p,p}_{xy}}
\norm{v'}_{X_{\beta_{z}}}, \\
\norm{F_z(Z,Z')}_{X^{-1}_{\overline{\sigma},p}}
&\leq C \norm{Z}_{H^{-1+1/2p,p}_zH^{1+1/p,p}_{xy}}
\norm{Z'}_{H^{1/2p,p}_z H^{1/p,p}_{xy}}.
\end{align*} 
\end{lemma}

\begin{proof}
Due to the boundedness of $\IP$, the definition of $H^{-1,p}_zL^p_{xy}$ which implies that $\partial_z\colon  L^{p}_zL^p_{xy} \rightarrow H^{-1,p}_zL^p_{xy}$ is a bounded operator, then by H\"older's inequality where 
$q_1, q_2\in (1,\infty)$ and $p_1, p_2\in (1,\infty)$ are such that $1/q_1+1/q_2=1/p$ and $1/p_1+1/p_2=1/p$ 
we estimate
with some $C>0$  
\begin{align*}
\norm{F_z(v,v')}_{X^{-1}_{\overline{\sigma},p}}
\leq C \norm{\partial_z (w(v) v')}_{H^{-1,p}_z L^p_{xy}}
\leq C \norm{ w(v) v'}_{L^p_z L^{p}_{xy}}  
\leq C \norm{w(v)}_{L^{p_1}_z L^{q_1}_{xy}} \norm{v'}_{L^{p_2}_z L^{q_2}_{xy}}.
\end{align*}
Next, we employ the notations $\int_0^z\phi$ for the function $(x,y,z)\mapsto \int_0^z\phi(x,y,\xi) d\xi$ and $\phi_0:=(\int_0^z\phi) - \overline{\phi}$ which satisfies $\partial_z \phi_0= \phi$ by the fundamental theorem of calculus with respect to the $z$-coordinate.
Then  by duality and using that $C_0^{\infty}(\cD)\subset L^{p_1'}_z L^{q_1'}_{xy}$ is dense, where $p_1'$ and $q_1'$ are such that  $1/p_1+1/p_1'=1$ and  $1/q_1+1/q_1'=1$, one obtains
\begin{align*}
\norm{w(v)}_{L^{p_1}_z L^{q_1}_{xy}} 
&= \sup\left\{ \abs{\int_{\cD}  w(v)\cdot   \phi} \colon \phi\in C_0^{\infty}(\cD) \hbox{ with } \norm{\phi}_{L^{p_1'}_z L^{q_1'}_{xy}}\leq 1 \right\}  
\\
&=\sup\left\{ \abs{\int_G\int_{-h}^0 w(v)\cdot   \partial_z\phi_0} \colon \phi\in C_0^{\infty}(\cD),  \norm{\partial_z\phi_0}_{L^{p_1'}_z L^{q_1'}_{xy}}\leq 1 \right\} \\ 
&= \sup\left\{ \abs{\int_G\left(\int_{-h}^0 \div_H v \cdot  \phi_0 + \left[ w(v)\cdot\phi_0\right]_{-h}^0\right)} \colon   \phi\in C_0^{\infty}(\cD),  \norm{\partial_z\phi_0}_{L^{p_1'}_z L^{q_1'}_{xy}}\leq 1 \right\} \\ 
&\leq \sup\left\{ \abs{\int_{\cD} \div_H v \cdot  \psi } \colon \psi\in H^{1,p_1'}_z L^{q_1'}_{xy} \hbox{ with } \norm{\psi}_{H^{1,p_1'}_z L^{q_1'}_{xy}}\leq 1+C_P \right\} \\
&=(1+C_P)^{-1} \norm{\div_H v}_{H^{-1,p_1}_zL_{xy}^{q_1}}\\
&\leq (1+C_P)^{-1}\norm{v}_{H^{-1,p_1}_zH^{1,q_1}_{xy}}.
\end{align*}
Here, we have used an integration by parts together with $\partial_z w(v)=-\div_H v$ and the fact that the boundary term vanishes due to the boundary conditions $w(v)(\cdot,\cdot,-h)=w(v)(\cdot,\cdot,0)=0$. Moreover it has been used that  $\phi_0\in H^{1,p_1'}_z L^{q_1'}_{xy}$ where $\phi_0(\cdot,\cdot,-h)=\phi_0(\cdot,\cdot,0)=0$ and therefore $\phi_0$ satisfies a Poincar\'{e} inequality with constant $C_P>0$. 
In particular for $p_1=p_2=2p$ and $q_1=q_2=2p$, one obtains by the embeddings
\begin{align*}
H^{-1+1/2p,p}_z \hookrightarrow H^{-1,2p}_z , \quad 
H^{1/2p,p}_z \hookrightarrow L^{2p}_z, \quad \hbox{and} \quad 
H^{1/p,p}_{xy} \hookrightarrow L^{2p}_{xy},
\end{align*}
where the respective Sobolev indices 
agree, that
\begin{align*}
\norm{F_z(v,v')}_{X^{-1}_{\overline{\sigma},p}} &\leq C \norm{v}_{H^{-1+1/2p,p}_zH^{1+1/p,p}_{xy}} \norm{v'}_{H^{1/2p,p}_z H^{1/p,p}_{xy}}.
\end{align*}

Next, using the the mapping properties of $\Delta_z$ and $\Delta_H$ (and that these operators commute with $\IP$), one can express the anisotropic Sobolev norms by the following graph norms 
\begin{align*}
\norm{v}_{H^{-1+1/2p,p}_zH^{1+1/p,p}_{xy}}&=\norm{(- \Delta_z+1)^{1/4p}(- \Delta_H+1)^{1/2+1/2p} v}_{X^{-1}_{\overline{\sigma},p}}
\quad \hbox{and} \\ 
\norm{v'}_{H^{1/2p,p}_z H^{1/2p}_{xy}}
&=\norm{(- \Delta_z+1)^{1/2+1/4p}(- \Delta_H+1)^{1/2p} v}_{X^{-1}_{\overline{\sigma},p}}.
\end{align*}
By the mixed derivative theorem, see e.g. \cite[Corollary 4.5.10]{PruessSimonett}, which applies to $\Delta_z$ and $\Delta_H$, one gets that
\begin{align*}
D\left((- \Delta_z+1)^{\beta_{z}} \cap(- \Delta_H+1)^{\beta_{z}}\right)
\hookrightarrow  D\left((- \Delta_z+1)^{1/4p}(- \Delta_H+1)^{1/2+1/2p}\right)
\end{align*}
if for some $\beta_{z}>0$ and $\theta\in (0,1)$ one has $\theta\beta_{z}=1/4p$ and $(1-\theta)\beta_{z}=1/2+1/2p$. For instance this holds for
\begin{align*}
\beta_{z}=1/2+ 3/4p\quad \hbox{and}  \quad\theta=(1/4p)/\beta_{z}.
\end{align*}
Note that in $H^{-1,p}_zL^p_{xy}$ for $\beta_z\geq 0$
\begin{align}\label{eq:powereta}
D\left((- \Delta_z+1)^{\beta_z}\cap(- \Delta_H+1)^{\beta_z}\right)
&= H^{-1+2\beta_z,p}_zL^p_{xy} \cap H^{-1,p}_zH^{2\beta_z,p}_{xy} \\
&=
D\left((\Delta_z + \Delta_H)^{\beta_z}\right) \notag
\end{align}
since the domains of the fractional powers are given by complex interpolation and by \cite[Lemma 9.5]{EscherPruessSimonett2003}
or \cite[Theorem 3.1]{GGS1993}. 
Similarly to the above
\begin{align*}
D\left((- \Delta_z+1)^{\beta_{z}}\cap(- \Delta_H+1)^{\beta_{z}}\right)
\hookrightarrow  D\left((- \Delta_z+1)^{1/2+1/4p}(- \Delta_H+1)^{1/2p}\right)
\end{align*}
if for some $\theta\in (0,1)$ one has $(1-\theta)\beta_{z}=1/2+1/4p$ and $\theta\beta_{z}=1/2p$ which
holds for
\begin{align*}
\beta_{z}=1/2+ 3/4p \quad \hbox{and}  \quad\theta=(1/2p)/\beta_{z}.
\end{align*}
Hence using again \eqref{eq:powereta}, we obtain the first inequality
\begin{align*}
\norm{F_z(v,v')}_{X^{-1}_{\overline{\sigma},p}} &\leq C \norm{v}_{D((A_p^{(-1)})^{1/2+3/4p})} \norm{v'}_{D((A_p^{(-1)})^{1/2+3/4p})}.
\end{align*}
Replacing $v'$ by $Z'$ and/or $v$ by $Z$, the other estimates follow analogously. 
\end{proof}

\begin{lemma}[Estimate on the non-linearity part II]\label{lemma:nonlin_II}
Let $p\in(2,\infty)$, $s\in [-1+1/p,0]$, and $\varepsilon>0$, then there exists a constant $C>0$ depending only on $h,p$ such that for $\beta_{H,1}=1/2+|s|/2+1/p +\varepsilon/2$ and $\beta_{H,2}=1-|s|/2$
\begin{align*}
\norm{F_H(v,v')}_{X^{-1}_{\overline{\sigma},p}} &\leq C 
(\norm{v}_{X_{\beta_{H,1}}}\norm{v'}_{X_{\beta_{H,2}}} + \norm{v'}_{X_{\beta_{H,1}}}\norm{v}_{X_{\beta_{H,2}}}), \\
\norm{F_H(v,Z')}_{X^{-1}_{\overline{\sigma},p}} 
&\leq C 
(\norm{v}_{X_{\beta_{H,1}}}\norm{Z'}_{H^{s,p}_zH^{1,p}_{xy}} + \norm{Z'}_{H^{|s|,p}_zH^{2/p+\varepsilon,p}_{xy}}\norm{v}_{X_{\beta_{H,2}}}), \\ 
\norm{F_H(Z,v')}_{X^{-1}_{\overline{\sigma},p}} &\leq C  
(\norm{v'}_{X_{\beta_{H,1}}}\norm{Z}_{H^{s,p}_zH^{1,p}_{xy}} + \norm{Z}_{H^{|s|,p}_zH^{2/p+\varepsilon,p}_{xy}}\norm{v'}_{X_{\beta_{H,2}}}),\\ 
\norm{F_H(Z,Z')}_{X^{-1}_{\overline{\sigma},p}} &\leq C  
(\norm{Z'}_{H^{|s|,p}_zH^{2/p+\varepsilon,p}_{xy}}\norm{Z}_{H^{s,p}_zH^{1,p}_{xy}} + \norm{Z}_{H^{|s|,p}_zH^{2/p+\varepsilon,p}_{xy}}\norm{Z'}_{H^{s,p}_zH^{1,p}_{xy}}).
\end{align*} 
\end{lemma}

\begin{proof}
Note that $F_H(v,v')=\IP(v\nabla_H v' + v'\nabla_H v)$.
By the embedding
\begin{align*}
\IP H^{s,r}_zL^{p}_{xy}\hookrightarrow X^{-1}_{\overline{\sigma},p} \quad \hbox{for }s\geq -1+1/r-1/p \quad \hbox{and} \quad r \in (1,p],
\end{align*} 
and the boundedness of $\IP$, one obtains for some $C>0$ that
\begin{align*}
\norm{F_H(v,v')}_{X^{-1}_{\overline{\sigma},p}} \leq C 
\norm{F_H(v,v')}_{H^{s,r}_zL^{p}_{xy}}.
\end{align*} 
By the mixed derivative theorem $$D(A_p^{(-1)})\hookrightarrow D((1-\Delta_z)^{1/2}(1-\Delta_H)^{1/2}) =L^p_zH^{1,p}_{xy}.$$ For $p\geq2$ therefore  $\IP(v\nabla_H v' + v'\div_H v)$
is a regular distribution. Hence
the pairing $\langle\cdot,\cdot\rangle$ of $H^{-s,r'}_zL^{p'}_{xy}$ and $H^{s,r}_zL^{p}_{xy}$ for $s\leq 0$, where $1/p+1/p'=1$ and $1/r+1/r'=1$, is given as follows
\begin{align*}
\norm{F_H(v,v')}_{X^{-1}_{\overline{\sigma},p}} &=
\sup\{  \abs{\int_{\cD}  \phi \cdot  \IP(v\nabla_H v' + v'\div_H v)}\colon \phi\in \IP H^{-s,r'}_zL^{p'}_{xy}, \norm{\phi}_{H^{-s,r'}_zL^{p'}_{xy}}\leq 1 \}
\\
&=
\sup\{  \abs{\int \IP\phi \cdot  (v\nabla_H v' + v'\div_H v) }\colon \phi\in \IP H^{-s,r'}_zL^{p'}_{xy}, \norm{\phi}_{H^{-s,r'}_zL^{p'}_{xy}}\leq 1 \}\\
&=
\sup\{  \abs{\int  (\phi v) \cdot  \nabla_H v' +  (\phi v')\cdot \div_H v}\colon \phi\in \IP H^{-s,r'}_zL^{p'}_{xy}, \norm{\phi}_{H^{-s,r'}_zL^{p'}_{xy}}\leq 1 \}\\
&=
\sup\{  \bigg\vert\int  (-\Delta_z+1)^{-s/2}(\phi v)\cdot (-\Delta_z+1)^{s/2}\nabla_H v'   \\
&\quad+ (-\Delta_z+1)^{-s/2}(\phi v')\cdot (-\Delta_z+1)^{s/2}\div_H v\bigg\vert\colon \phi\in \IP H^{-s,r'}_zL^{p'}_{xy}, \norm{\phi}_{H^{-s,r'}_zL^{p'}_{xy}}\leq 1 \}\\
&\leq \sup\{\norm{\varphi}_{H^{-s,r'}_zL_{xy}^{p'}}\colon \phi\in \IP H^{-s,r'}_zL^{p'}_{xy}, \norm{\phi}_{H^{-s,r'}_zL^{p'}_{xy}}\leq 1 \} \\ &\qquad\cdot(
\norm{v}_{H^{-s,r_1}_zL_{xy}^{p_1}}
\norm{\nabla_H v'}_{H^{s,r_2}_zL_{xy}^{p_2}} + 
\norm{v'}_{H^{-s,r_1}_zL_{xy}^{p_1}}
\norm{\nabla_H v}_{H^{s,r_2}_zL^{p_2}_{xy}}),	
\end{align*}
where $1/r'+1/r_1+1/r_2=1$. Here we have used that that by duality for $f$ and $g$ sufficiently regular
\begin{align*}
\int f\cdot g =\int  (-\Delta_z+1)^{s/2}(-\Delta_z+1)^{-s/2}f\cdot g 
=\int  (-\Delta_z+1)^{-s/2}f\cdot (-\Delta_z+1)^{s/2} g.
\end{align*}
Taking $r=p/2$ for $p>2$ and $r_1=r_2=p$, $p_1=\infty$ and $p_2=p$, we obtain by the  Sobolev embedding $H^{2/p+\varepsilon,p}_{xy}\hookrightarrow L^\infty_{xy}$  for any $\varepsilon>0$
that 
\begin{align*}
\norm{F_H(v,v')}_{X^{-1}_{\overline{\sigma},p}} &\leq
C(  
\norm{v}_{H^{|s|,p}_zL^{\infty}_{xy}}
\norm{\nabla_H v'}_{H^{s,p}_zL^{p}_{xy}}
+
\norm{v'}_{H^{|s|,p}_zL^{\infty}_{xy}}
\norm{\div_H v}_{H^{s,p}_zL^{p}_{xy}}) \\
&\leq C \left(
\norm{v}_{H^{|s|,p}_zH^{2/p+\varepsilon,p}_{xy}}
\norm{ v'}_{H^{s,p}_zH^{1,p}_{xy}}
+
\norm{v'}_{H^{|s|,p}_zH^{2/p+\varepsilon,p}_{xy}}
\norm{v}_{H^{s,p}_zH^{1,p}_{xy}}\right),
\end{align*}
where we write $|s|=-s$ to clarify the signs.
As before, one can express then
\begin{align*}
\norm{v}_{H^{-s,p}_zH^{2/p+\varepsilon,p}_{xy}}&=\norm{(- \Delta_z^{(-1)}+1)^{1/2+|s|/2}(- \Delta_H+1)^{1/p +\varepsilon/2} v}_{X^{-1}_{\overline{\sigma},p}}
\quad \hbox{and}, \\ 
\norm{v}_{H^{s,p}_zH^{1,p}_{xy}}
&=\norm{(- \Delta_z+1)^{1/2-|s|/2}(- \Delta_H+1)^{1/2}}_{X^{-1}_{\overline{\sigma},p}}.
\end{align*}
By the mixed derivative theorem
\begin{align*}
D\left((- \Delta_z+1)^{\beta_{H,1}}\cap(- \Delta_H+1)^{\beta_{H,1}}\right)
\hookrightarrow  D\left((- \Delta_z^{(-1)}+1)^{1/2+|s|/2}(- \Delta_H+1)^{1/p +\varepsilon/2}\right)
\end{align*}
if for some $\theta\in (0,1)$ one has $\theta\beta_{H,1}=1/2+|s|/2$ and $(1-\theta)\beta_{H,1}=1/p +\varepsilon/2$, which holds for
\begin{align*}
\beta_{H,1}=1/2+|s|/2+1/p +\varepsilon/2\quad \hbox{and}  \quad\theta=(1/2+|s|/2)/\beta_{H,1}.
\end{align*}
Furthermore,
\begin{align*}
D\left((- \Delta_z+1)^{\beta_{H,2}}\cap(- \Delta_H+1)^{\beta_{H,2}}\right)
\hookrightarrow (- \Delta_z+1)^{1/2-|s|/2}(- \Delta_H+1)^{1/2}
\end{align*}
if for some $\theta\in (0,1)$ one has $\theta\beta_{H,2}=1/2-|s|/2$ and $(1-\theta)\beta_{H,2}=1/2$, which holds for
\begin{align*}
\beta_{H,2}=1-|s|/2\quad \hbox{and}  \quad\theta=(1/2-|s|)/\beta_{H,2}.
\end{align*}
Hence,  we can estimate
\begin{align*}
\norm{F_H(v,v')}_{X^{-1}_{\overline{\sigma},p}} \leq C 
(\norm{v}_{D((A_p^{(-1)})^{\beta_{H,1}})}\norm{v'}_{D((A_p^{(-1)})^{\beta_{H,2}})} + \norm{v'}_{D((A_p^{(-1)})^{\beta_{H,1}})}\norm{v}_{D((A_p^{(-1)})^{\beta_{H,2}})}).
\end{align*}
This proves the first estimate and combining it with the previous estimates, the other estimates follow analogously. 
\end{proof}

\subsection{Local well-posedness}
\begin{proof}[Proof of Theorem~\ref{thm:main1}]
Here, we apply the result from \cite{PruessWilke2016} to solve ~\eqref{eq:PE_stoch_z}. Due to Proposition~\ref{prop:Hinfty} the hydrostatic Stokes operator has a bounded $H^\infty$-calculus and hence to apply result from \cite{PruessWilke2016}, we have only to add suitable estimates on the non-linearity. 
By the bilinearity of $F(\cdot,\cdot)$ one has 
\begin{multline*}
F(Z+v,Z+v)-F(Z+v',Z+v')=  F(v-v',v) + F(v',v-v') + F(Z,v-v') + F(v-v',Z).
\end{multline*}
To control the linear terms $F(Z,v-v') + F(v-v',Z)$ and the force term $F(Z,Z)$, we can use first the estimates from Lemma~\ref{lemma:nonlin_I} and Lemma~\ref{lemma:nonlin_II}.
Then
the regularity results from Proposition~\ref{prop:reg_Zb} allow us to assure the necessary regularity of $Z$.  
By shifting the scales also in the horizontal direction,
\begin{align*}
\norm{Z}_{H^{|s|,p}_zH^{2/p+\varepsilon,p}_{xy}}&=\norm{(- \Delta_z^{(-1)}+1)^{1/2} v}_{X^{-1+|s|, 2/p+\varepsilon}_{\overline{\sigma},p}}
\quad \hbox{and}, \\ 
\norm{Z}_{H^{s,p}_zH^{1,p}_{xy}}
&=\norm{(- \Delta_z+1)^{1/2-|s|}(- \Delta_H+1)^{1/2-1/p-\varepsilon/2}}_{X^{-1+|s|, 2/p+\varepsilon}_{\overline{\sigma},p}},\\
\norm{Z}_{H^{-1+1/2p,p}_zH^{1+1/p,p}_{xy}}&=\norm{(- \Delta_z+1)^{\max\{1/4p-|s|/2,0\}}(- \Delta_H+1)^{1/2-1/2p-\varepsilon/2}}_{X^{-1+|s|, 2/p+\varepsilon}_{\overline{\sigma},p}},  \\
\norm{Z'}_{H^{1/2p,p}_z H^{1/p,p}_{xy}}&=\norm{(- \Delta_z+1)^{1/2+1/4p-|s|/2}}_{X^{-1+|s|, 2/p+\varepsilon}_{\overline{\sigma},p}}.
\end{align*}

By the mixed derivative theorem
\begin{align*}
D\left((-\Delta_z+1)^{1/2}\cap(- \Delta_H+1)^{1/2}\right)
\hookrightarrow  D\left((- \Delta_z^{(-1)}+1)^{1/2-|s|/2}(- \Delta_H+1)^{1/2-1/p -\varepsilon/2}\right)
\end{align*}
if for some $\theta\in (0,1)$ one has $\theta/2=1/2-|s|$ and $(1-\theta)/2=1/2-1/p -\varepsilon/2$, which holds for
\begin{align*}
\theta/2=1/2-|s| 
\quad \hbox{and}  \quad(1-\theta)/2= 1/2-1/p -\varepsilon/2 
\end{align*}
if we chose $|s|= 1/2 -1/p-\varepsilon/2$.
Moreover, we require by Proposition~\ref{prop:reg_Zb} that
\begin{align}\label{eq:pleq4}
|s|= 1/2 -1/p-\varepsilon/2 <1/p  \quad \hbox{that is}\quad p\leq 4,
\end{align} 
where $\varepsilon$ has been chosen sufficiently small. Moreover, 
with this choice
\begin{align*}
D\left((-\Delta_z+1)^{1/2}\cap(- \Delta_H+1)^{1/2}\right)
\hookrightarrow        
D((- \Delta_z+1)^{\max\{1/4p-|s|/2,0\}}(- \Delta_H+1)^{1/2-1/2p-\varepsilon/2}),
\end{align*}
then
\begin{align}\label{eq:pge3}
D\left((-\Delta_z+1)^{1/2}\cap(- \Delta_H+1)^{1/2}\right)
\hookrightarrow        
D((- \Delta_z+1)^{1/2+1/4p-|s|/2}) \quad \hbox{for }p>3.
\end{align}

Hence, we can estimate using Lemma~\ref{lemma:nonlin_I} and Lemma~\ref{lemma:nonlin_II} to obtain with $\beta=\max\{\beta_z,\beta_{H,1},\beta_{H,2}\}$
\begin{multline*}
\norm{F(Z+v,Z+v)-F(Z+v',Z+v')}_{X^{-1}_{\overline{\sigma},p}} \leq C (\norm{v}_{X_\beta} + \norm{v'}_{X_\beta})\norm{v-v'}_{X_\beta} \\
+ C\norm{Z}_{D((A_p^{(-1+|s|, 2/p+\varepsilon)})^{1/2})} \norm{v-v'}_{X_\beta}
=\sum_{j=1}^2
C_j(\norm{v}^{\rho_j}_{X_\beta} + \norm{v'}^{\rho_j}_{X_\beta})\norm{v-v'}_{X_{\beta}},
\end{multline*}
where 
$\rho_1=1$ while $\rho_2=0$, and $C_1=C$ and $C_2= C\norm{Z}_{D((A_p^{(-1+|s|, 2/p+\varepsilon)})^{1/2})}$.
In proof of the result from \cite{PruessWilke2016}, one just has to shorten the estimate on the terms involving $Z$ to become
\begin{align*}
\left(\int_0^T\norm{Z}_{D((A_p^{(-1+|s|, 2/p+\varepsilon)})^{1/2})} \norm{v-v'}_{X_\beta} t^{(1-\mu)q}\right)^{1/q} \\
\leq \norm{Z}_{L_{\sigma'}^{2q}(0,T; D((A_p^{(-1+|s|, 2/p+\varepsilon)})^{1/2}))}  \norm{v-v'}_{L_{\sigma}^{2q}(0,T;X_{\beta})}, 
\end{align*}
where $\sigma=\sigma'=(1+\mu)/2$.
Then one proceeds as there by the estimate
\begin{equation*} 
\norm{v-v'}_{L_{\sigma}^{2q}(0,T;X_{\beta})} \leq C\norm{v-v'}_{H_{\mu}^{1-\beta, q}(0,T;X_{\beta})}
\end{equation*} 
while keeping the term $\norm{Z}_{L_{\sigma'}^{2q}(0,T; D((A_p^{(-1+|s|, 2/p+\varepsilon)})^{1/2})}$. Note that with the choice of $s$ and $p$ one has $\beta_z<\beta=\beta_{H,1}=\beta_{H,2}$.
Hence the critical weight is 
\begin{align*}
\mu_c=1/q+3\beta/2-1/2 \quad \hbox{with} \quad \beta=3/4+1/2p+\varepsilon/4
\end{align*}
as desired.
\end{proof}

\subsection*{Acknowledgements}
Matthias Hieber gratefully acknowledges the support by the Deutsche Forschungsgemeinschaft (DFG) through  the Research Unit 5528. Amru Hussein has been generously supported by 
Deutsche Forschungsgemeinschaft (DFG) -- project number 508634462 and by MathApp -- Mathematics Applied to Real-World Problems -- part of the Research Initiative of the Federal State of Rhineland-Palatinate, Germany. Martin Saal  gratefully acknowledges the financial support of the Deutsche Forschungsgemeinschaft (DFG) through the 
research fellowship SA 3887/1-1.

\bibliographystyle{alpha}
\bibliography{literature.bib}

\begin{bibdiv}
\begin{biblist}

\bib{agresti2023global}{unpublished}{
      author={Agresti, A.},
      author={Luongo, E.},
       title={Global well-posedness and interior regularity of {2D}
  {N}avier-{S}tokes equations with stochastic boundary conditions},
        date={2023},
         url={arXiv:2306.11081},
}

\bib{AgrestiVeraar}{article}{
      author={Agresti, Antonio},
      author={Veraar, Mark},
       title={Nonlinear parabolic stochastic evolution equations in critical
  spaces part i. stochastic maximal regularity and local existence},
        date={2022},
     journal={Nonlinearity},
      volume={35},
      number={8},
       pages={4100},
}

\bib{AgrestiVeraarII}{article}{
      author={Agresti, Antonio},
      author={Veraar, Mark},
       title={Nonlinear parabolic stochastic evolution equations in critical
  spaces part ii.},
        date={2022},
     journal={J. Evol. Equ.},
      volume={22},
      number={56},
         url={https://dx.doi.org/10.1007/s00028-022-00786-7},
}

\bib{Amann1993}{incollection}{
      author={Amann, H.},
       title={Nonhomogeneous linear and quasilinear elliptic and parabolic
  boundary value problems},
        date={1993},
   booktitle={Function spaces, differential operators and nonlinear analysis},
      editor={Schmeisser, H.},
      editor={Triebel, H.},
   publisher={Vieweg+Teubner Verlag},
       pages={9\ndash 126},
}

\bib{Amann2019}{book}{
      author={Amann, H.},
       title={Linear and quasilinear parabolic problems. {V}ol. {II}},
   publisher={Birkh\"auser/Springer, [Cham]},
        date={2019},
}

\bib{ABHN}{book}{
      author={Arendt, W.},
      author={Batty, Ch.~J.~K.},
      author={Hieber, M.},
      author={Neubrander, F.},
       title={Vector-valued {L}aplace transforms and {C}auchy problems},
   publisher={Birkh\"{a}user/Springer Basel AG, Basel},
        date={2011},
}

\bib{AB:02}{article}{
      author={Arendt, W.},
      author={Bu, S.},
       title={The operator-valued marcinkiewicz multiplier theorem and maximal
  regularity},
        date={2002},
     journal={Math. Z.},
      volume={240},
      number={2},
       pages={311\ndash 343},
}

\bib{BE20}{article}{
      author={Binz, T.},
      author={Engel, K.-J.},
       title={First order evolution equations with dynamic boundary
  conditions},
        date={2020},
     journal={Phil. Trans. A},
}

\bib{BS01}{article}{
      author={Bresch, D.},
      author={Simon, J.},
       title={On the effect of friction on wind driven shallow lakes},
        date={2001},
     journal={J. Math. Fluid Mech.},
      volume={3},
       pages={231\ndash 258},
}

\bib{BS20}{article}{
      author={Brzeźniak, Zdzisław},
      author={Slavík, Jakub},
       title={Well-posedness of the 3d stochastic primitive equations with
  multiplicative and transport noise},
        date={2021},
        ISSN={0022-0396},
     journal={Journal of Differential Equations},
      volume={296},
       pages={617\ndash 676},
}

\bib{CaoTiti2007}{article}{
      author={Cao, Ch.},
      author={Titi, E.},
       title={Global well--posedness of the three-dimensional viscous primitive
  equations of large scale ocean and atmosphere dynamics},
        date={2007},
     journal={Annals of Mathematics},
      volume={166},
       pages={245\ndash 267},
}

\bib{CENN03}{article}{
      author={Casarino, V.},
      author={Engel, K.-J.},
      author={Nagel, R.},
      author={Nickel, G.},
       title={A semigroup approach to boundary feedback systems},
        date={2003},
     journal={Integral Equations Operator Theory},
      volume={47},
       pages={289\ndash 306},
}

\bib{CL01}{article}{
      author={Cessi, P.},
      author={Louazel, S.},
       title={Decadal oceanic response to stochastic wind forcing},
        date={2001},
     journal={J. Physical Oceanography},
      volume={34},
       pages={3020\ndash 3029},
}

\bib{DZ96}{book}{
      author={Da~Prato, G.},
      author={Zabcyyk, J.},
       title={Ergodicity for infinite dimensional systems},
      series={London Math. Soc. Lecture Note Series},
   publisher={Cambridge Univ. Press},
        date={1996},
      volume={229},
}

\bib{DS09}{article}{
      author={Dalibard, A.},
      author={Saint-Raymond, L.},
       title={Resonant wind-driven oceanic motions},
        date={2009},
     journal={C.R. Acad. Sci Paris},
      volume={347},
       pages={451\ndash 456},
}

\bib{DebusscheGlattholtzTemam2011}{article}{
      author={Debussche, A.},
      author={Glatt-Holtz, N.},
      author={Temam, R.},
       title={Local martingale and pathwise solutions for an abstract fluids
  model},
        date={2011},
     journal={Physica D: Nonlinear Phenomena},
      volume={240},
       pages={1123\ndash 1144},
}

\bib{DebusscheGlattholtzTemamZiane2012}{article}{
      author={Debussche, A.},
      author={Glatt-Holtz, N.},
      author={Temam, R.},
      author={Ziane, M.},
       title={Global existence and regularity for the 3d stochastic primitive
  equations of the ocean and atmosphere with multiplicative white noise},
        date={2012},
     journal={Nonlinearity},
      volume={25},
      number={7},
       pages={2093\ndash 2118},
}

\bib{DG00}{article}{
      author={Desjardin, B.},
      author={Grenier, E.},
       title={On the homogeneous model of wind-driven ocean circulation},
        date={2000},
     journal={SIAM J. Appl. Math.},
      volume={60},
       pages={43\ndash 60},
}

\bib{DongZhang2017}{article}{
      author={Dong, Z.},
      author={Zhang, R.},
       title={Markov selection and w-strong feller for 3d stochastic primitive
  equations},
        date={2017},
     journal={Sci. China Math.},
      volume={60},
       pages={1873\ndash 1900},
}

\bib{EscherPruessSimonett2003}{article}{
      author={Escher, J.},
      author={Pr\"{u}ss, J.},
      author={Simonett, G.},
       title={Analytic solutions for a {S}tefan problem with {G}ibbs-{T}homson
  correction},
        date={2003},
     journal={J. Math. Phys.},
      volume={563},
       pages={1\ndash 52},
}

\bib{GGS1993}{inproceedings}{
      author={Giga, Mariko},
      author={Giga, Yoshikazu},
      author={Sohr, Hermann},
       title={Lp estimates for the stokes system},
        date={1993},
   booktitle={Functional analysis and related topics, 1991},
      editor={Komatsu, Hikosaburo},
   publisher={Springer Berlin Heidelberg},
     address={Berlin, Heidelberg},
       pages={55\ndash 67},
}

\bib{GGHHK17}{article}{
      author={Giga, Y.},
      author={Gries, M.},
      author={Hieber, M.},
      author={Hussein, A.},
      author={Kashiwabara, T.},
       title={Bounded {$H^\infty$}-calculus for the hydrostatic stokes operator
  on {$L^p$}-spaces and applications},
        date={2017},
     journal={Proc. Amer. Math. Soc.},
      volume={145},
      number={9},
       pages={3865\ndash 3876},
}

\bib{GGHHK20}{article}{
      author={Giga, Y.},
      author={Gries, M.},
      author={Hieber, M.},
      author={Hussein, A.},
      author={Kashiwabara, T.},
       title={Analyticity of solutions of primitive equations},
        date={2020},
     journal={Math. Nachr.},
      volume={293},
      number={2},
       pages={284\ndash 304},
}

\bib{GGHHK20b}{article}{
      author={Giga, Y.},
      author={Gries, M.},
      author={Hieber, M.},
      author={Hussein, A.},
      author={Kashiwabara, T.},
       title={The hydrostatic stokes operator and well-posedness of the
  primitive equations on spaces of bounded functions},
        date={2020},
     journal={J. Funct. Anal.},
      volume={279},
      number={108561},
       pages={3865\ndash 3876},
}

\bib{BCP15}{article}{
      author={Giuseppe~Buffoni, Andrea~Cappelletti},
      author={Picco, Paola},
       title={On the ekman equations for ocean currents driven by a stochastic
  wind},
        date={2015},
     journal={Stochastic Analysis and Applications},
      volume={33},
      number={2},
       pages={356\ndash 382},
}

\bib{GlattholtzKukavicaVicolZiane2014}{article}{
      author={Glatt-Holtz, N.},
      author={Kukavica, I.},
      author={Vicol, V.},
      author={Ziane, M.},
       title={Existence and regularity of invariant measures for the three
  dimensional stochastic primitive equations},
        date={2014},
     journal={J. Math. Phys.},
      volume={55},
}

\bib{GlattholtzTemam2011}{article}{
      author={Glatt-Holtz, N.},
      author={Temam, R.},
       title={Pathwise solutions of the 2--d stochastic primitive equations},
        date={2011},
     journal={Appl. Math. Optim.},
      volume={63},
      number={3},
       pages={401\ndash 433},
}

\bib{Gre87}{article}{
      author={Greiner, G.},
       title={Perturbing the boundary conditions of a generator},
        date={1987},
     journal={Houston J. Math},
      volume={13},
       pages={213\ndash 229},
}

\bib{GM97}{article}{
      author={Grenier, E.},
      author={Masmoudi, N.},
       title={Ekman layers of rotating fluids: the case of well-prepared
  initial data},
        date={1997},
     journal={Comm. Partial Diff. Equ.},
      volume={22},
       pages={953\ndash 975},
}

\bib{GuoHuang2009}{article}{
      author={Guo, B.},
      author={Huang, D.},
       title={{3D} stochastic primitive equations of the large-scale ocean:
  global well-posedness and attractors},
        date={2009},
     journal={Comm. Math. Phys.},
      volume={286},
       pages={697\ndash 723},
}

\bib{HHK2016}{article}{
      author={Hieber, M.},
      author={Hussein, A.},
      author={Kashiwabara, T.},
       title={Global strong {$L^p$} well-posedness of the 3d primitive
  equations with heat and salinity diffusion},
        date={2016},
     journal={J. Differ. Equ.},
      volume={261},
      number={12},
       pages={6950\ndash 6981},
}

\bib{HK16}{article}{
      author={Hieber, M.},
      author={Kashiwabara, T.},
       title={Global strong well–posedness of the three dimensional primitive
  equations in {$L^p$}–spaces},
        date={2016},
     journal={Arch. Rational Mech. Anal.},
      volume={221},
      number={3},
       pages={1077\ndash 1115},
}

\bib{HH20}{incollection}{
      author={Hieber, Matthias},
      author={Hussein, Amru},
       title={An approach to the primitive equations for oceanic and
  atmospheric dynamics by evolution equations},
        date={2020},
   booktitle={Fluids under pressure},
      editor={Bodn{\'a}r, Tom{\'a}{\v{s}}},
      editor={Galdi, Giovanni~P.},
      editor={Ne{\v{c}}asov{\'a}, {\v{S}}{\'a}rka},
   publisher={Springer International Publishing},
     address={Cham},
       pages={1\ndash 109},
}

\bib{Hu2020}{article}{
      author={Hussein, A.},
       title={Partial and full hyper-viscosity for navier-stokes and primitive
  equations},
        date={2020},
     journal={J. Differ. Equ.},
      volume={269},
      number={4},
       pages={3003 \ndash  3030},
}

\bib{Hyt2017}{book}{
      author={Hyt\"{o}nen, T.},
      author={van Neerven, J.},
      author={Veraar, M.},
      author={Weis, L.},
       title={Analysis in {B}anach spaces. {V}ol. {II}},
   publisher={Springer, [Cham]},
        date={2017},
}

\bib{Ju2017}{article}{
      author={Ju, N.},
       title={On {$H^2$} solutions and {$z$}-weak solutions of the 3{D}
  primitive equations},
        date={2017},
     journal={Indiana Univ. Math. J.},
      volume={66},
      number={3},
       pages={973\ndash 996},
}

\bib{KalWei01}{article}{
      author={Kalton, N.},
      author={Weis, L.},
       title={The ${H}^\infty$-calculus and sums of closed operators},
        date={2001},
     journal={Math. Ann.},
      volume={321},
       pages={319\ndash 345},
}

\bib{KZ07}{article}{
      author={Kukavica, I.},
      author={Ziane, M.},
       title={On the regularity of the primitive equations of the ocean},
        date={2007},
     journal={Nonlinearity},
      volume={20},
       pages={2739\ndash 2753},
}

\bib{LiTiti2016}{incollection}{
      author={Li, J.},
      author={Titi, E.},
       title={Recent advances concerning certain class of geophysical flows},
        date={2016},
   booktitle={Handbook of mathematical analysis in mechanics of viscous
  fluids},
   publisher={Springer International Publishing},
}

\bib{Lionsetal1992}{article}{
      author={Lions, J.~L.},
      author={Temam, R.},
      author={Wang, Sh.~H.},
       title={New formulations of the primitive equations of atmosphere and
  applications},
        date={1992},
     journal={Nonlinearity},
      volume={5},
       pages={237\ndash 288},
}

\bib{Lionsetal1992_b}{article}{
      author={Lions, J.~L.},
      author={Temam, R.},
      author={Wang, Sh.~H.},
       title={On the equations of the large-scale ocean},
        date={1992},
     journal={Nonlinearity},
      volume={5},
       pages={1007\ndash 1053},
}

\bib{Lionsetal1993}{article}{
      author={Lions, J.~L.},
      author={Temam, R.},
      author={Wang, Sh.~H.},
       title={Models for the coupled atmosphere and ocean. ({CAO} {I},{II})},
        date={1993},
     journal={Comput. Mech. Adv.},
      volume={1},
       pages={3\ndash 119},
}

\bib{Nau2012}{thesis}{
      author={Nau, T.},
       title={{$L^p$}-theory of cylindrical boundary value problems. an
  operator-valued fourier multiplier and functional calculus approach.},
        type={Ph.D. Thesis},
        date={2012},
}

\bib{PTZ09}{incollection}{
      author={Petcu, M.},
      author={Temam, R.},
      author={Ziane, M.},
       title={Some mathematical problems in geophysical fluid dynamics},
        date={2008},
   booktitle={Handbook of numerical analysis},
      editor={Ciarlet, P.G.},
   publisher={Elsevier},
}

\bib{PruessSimonett}{book}{
      author={Pr\"uss, J.},
      author={Simonett, G.},
       title={Moving interfaces and quasilinear parabolic evolution equations},
   publisher={Birkh\"auser/Springer, [Cham]},
        date={2016},
}

\bib{PruessSimonettWilke}{article}{
      author={Pr{\"u}ss, J.},
      author={Simonett, G.},
      author={Wilke, M.},
       title={Critical spaces for quasilinear parabolic evolution equations and
  applications},
        date={2018},
     journal={J. Differ. Equ.},
      volume={264},
       pages={2028\ndash 2074},
}

\bib{PruessWilke2016}{article}{
      author={Pr{\"u}ss, J.},
      author={Wilke, M.},
       title={Addendum to the paper "{O}n quasilinear parabolic evolution
  equations in weighted {$L^p$}-spaces {II}"},
        date={2017},
     journal={J. Evol. Eq.},
      volume={7},
      number={4},
       pages={1381\ndash 1388},
}

\bib{SchnaubeltVeraar2011}{incollection}{
      author={Schnaubelt, Roland},
      author={Veraar, Mark},
       title={Stochastic equations with boundary noise},
        date={2011},
   booktitle={Parabolic problems: The herbert amann festschrift},
      editor={Escher, Joachim},
      editor={Guidotti, Patrick},
      editor={Hieber, Matthias},
      editor={Mucha, Piotr},
      editor={Pr{\"u}ss, Jan~W.},
      editor={Shibata, Yoshihiro},
      editor={Simonett, Gieri},
      editor={Walker, Christoph},
      editor={Zajaczkowski, Wojciech},
   publisher={Springer Basel},
     address={Basel},
       pages={609\ndash 629},
}

\bib{Ste93}{book}{
      author={Stein, E.M.},
       title={Harmonic analysis},
   publisher={Princeton University Press, Princeton},
        date={1993},
}

\bib{TZ96}{article}{
      author={Temam, R.},
      author={Ziane, M.},
       title={Navier-stokes equations in thin domains with various boundary
  conditions},
        date={1996},
     journal={Adv. Diff. Equ.},
      volume={1},
       pages={499\ndash 546},
}

\bib{Triebel}{book}{
      author={Triebel, H.},
       title={Theory of function spaces (reprint of 1983 edition)},
   publisher={Springer AG, Basel},
        date={2010},
}

\bib{NeervenVeraarWeis2}{article}{
      author={van Neerven, Jan},
      author={Veraar, Mark},
      author={Weis, Lutz},
       title={Maximal {$L^p$}-regularity for stochastic evolution equations},
        date={2012},
     journal={SIAM Journal on Mathematical Analysis},
      volume={44},
      number={3},
       pages={1372\ndash 1414},
}

\bib{NeervenVeraarWeis}{article}{
      author={van Neerven, Jan},
      author={Veraar, Mark},
      author={Weis, Lutz},
       title={Stochastic maximal $l^p$-regularity},
        date={2012},
     journal={Ann. Probab.},
      volume={40},
       pages={788\ndash 812},
}

\end{biblist}
\end{bibdiv}

\end{document}